\documentclass[xypic,amscd,syntonly,amssymb]{amsart}

\usepackage{graphicx}
\usepackage[all]{xy}
\usepackage{amssymb}
\usepackage{amsmath}

\usepackage{epsfig}
\theoremstyle{plain}
\newtheorem{Thm}{Theorem}
\newtheorem{Prop}[Thm]{Proposition}
\newtheorem{Cor}[Thm]{Corollary}
\newtheorem{Lem}[Thm]{Lemma}

 \theoremstyle{definition}
\theoremstyle{remark}

\numberwithin{equation}{section}

\begin{document}
 \title{Continuous families of  Hamiltonian torus actions}

 \author{ ANDR\'{E}S   VI\~{N}A}
\address{Departamento de F\'{i}sica. Universidad de Oviedo.   Avda Calvo
 Sotelo.     33007 Oviedo. Spain. }
 \email{vina@uniovi.es}
\thanks{This work has been partially supported by Ministerio de Ciencia y
Tecnolog\'{\i}a, grant MAT2003-09243-C02-01}

 \maketitle
\begin{abstract}
 We determine  conditions under which two Hamiltonian torus actions
 on a symplectic manifold $M$ are homotopic by a family of
 Hamiltonian torus actions, when  $M$ is a
 toric manifold  and when  $M$ is a coadjoint orbit.

\end{abstract}
\section {Introduction} \label{S:intro}

Let $(M,\omega)$ be a closed symplectic $2n$-manifold, and let $G$
be a compact Lie group. Given two Hamiltonian $G$-actions
on $M$, it is a natural issue to analyze  if one can be deformed
in the another by a continuous family of Hamiltonian $G$-actions.
In this   paper we will consider this matter; that is,  we will
 determine necessary conditions for two Hamiltonian
$G$-actions on $M$ to be connected by a homotopy consisting of
Hamiltonian $G$-actions.

The mathematical setting is the following:
  ${\rm
Ham}(M)$  will  be  the Hamiltonian group of $(M,\omega)$
\cite{Mc-S}, \cite{lP01}, and ${\rm Hom}(G,\,{\rm Ham}(M))$ will
denote the space of   all Lie group homomorphisms $\psi$ from $G$
to ${\rm Ham}(M)$; that is, the space consisting of the
Hamiltonian $G$-actions on $M$.

 Given $\psi,\psi'\in{\rm Hom}(G,\, {\rm Ham}(M))$, we write
$\psi\sim \psi'$ iff $\psi$ and $\psi'$ belong to the same
connected component of $\text {Hom} (G,\, \text {Ham}(M))$. In
other words, $\psi\sim \psi'$ iff there is a continuous family
$\{\psi^s: G\to {\rm  Ham }(M)) \}_{s\in[0,1]}$ of Lie group
homomorphisms, such that $\psi^0=\psi$ and $\psi^1=\psi'$; that
is, $\psi$ and $ \psi'$ are homotopic by a family of {\em group
homomorphisms}. The corresponding quotient space is denoted by
$[G,\,\text {Ham}(M)]_{gh}$.

In this paper we will consider the following three general
aspects:

  {\it 1. Necessary conditions for the $\sim$-equivalence.}
   It would be desirable to obtain necessary
and sufficient conditions for the $\sim$-equivalence of two
$G$-actions in the general case, but we will prove  more modest
results.    To each Hamiltonian $G$-action $\psi$ we will
associate a cohomology class $S_{\psi}$, such that the map $S_{-}$
is compatible with the relation $\sim$, and we will use the map
$S_{-}$  to distinguish different elements of $[G,\,\text
{Ham}(M)]_{gh}$ (see Theorem \ref{ThmHHA1}). In some particular
cases we will give direct proofs of the $\sim$-equivalence of two
actions on which the map $S_{-}$ takes the same  value.

\smallskip

{\it 2. Equivalence relations between toric actions.}
  A $G$-action can be regarded  as a
 representation of $G$ in the group ${\rm Ham}(M)$, so we can
 consider another natural  equivalence relation in the space ${\rm Hom}(G,\, {\rm Ham}(M))$:
 Two Hamiltonian $G$-actions on $M$,
 $\psi$ and $\psi'$, are said to be {\em  $r$-equivalent} if there exist
an element $h\in{\rm Ham}(M)$, such that $h\circ\psi_g\circ
h^{-1}=\psi'_g$ for all $g\in G$; that is, if they are equivalent
as representations of $G$ in the group $\text {Ham} (M)$. We will
compare the $r$-equivalence and the relation $\sim$ in the space
of   toric actions on $M$.  We  will prove, among other
properties, that the $r$-equivalence implies the
$\sim$-equivalence, and  that the equality of the moment polytopes
of $\psi$ and $\psi'$ is a necessary condition for the
$\sim$-equivalence   of them (see Theorem \ref{FinalThm}).

\smallskip

 {\it 3. $\sim$-equivalence of $U(1)$-actions.}
  When $G=U(1)$,  each $\psi\in{\rm Hom}(G,\, {\rm Ham}(M))$
 determines an element of $\{\psi\}\in\pi_1({\rm Ham}(M))$.
 We will compare our relation $\sim$ with the equivalence
 in $\pi_1({\rm Ham}(M))$ in some particular cases for which the homotopy type of ${\rm Ham}(M)$ is known.
 In \cite{Gr} Gromov proved that ${\rm Ham }({\mathbb
 C}P^2,\,\omega_{FS})$ is homotopy equivalent to the group
 $PU(3)$, where $\omega_{FS}$ is the Fubini-Study  symplectic
 structure, and  in \cite{Ab-McD} Abreu and McDuff determined the
 homology of the group of   symplectomorphisms of a rational ruled
 surface. In order to carry out the aforementioned comparison
 and also by its very interest,
we will   study  the relation $\sim$ between $U(1)$-actions
in two cases: \\
\hspace*{20pt}a) When   $M$ is the total space of the fibration
\begin{equation}\label{DefPL}
{\mathbb P}(L_1\oplus\dots\oplus L_{n-1}\oplus{\mathbb
C})\to{\mathbb
 C}P^1,
  \end{equation}
\hspace*{29pt} $L_i$ being a holomorphic line bundle over
${\mathbb C}P^1$. \\
 \hspace*{20pt}b) When the manifold $M$ is a coadjoint orbit; we look more
closely at the case of $U(1)$-actions on coadjoint orbits of
$SU(n)$ diffeomorphic to Grassmann manifolds.

The comparison of  the
 results of Gromov and Abreu-McDuff  above
 mentioned
  and our results   shows
 the  non-injectivity of the obvious map
  \begin{equation}\label{injective}
   [\psi]\in[U(1),\,\text
 {Ham}(M)]_{gh}\mapsto \{\psi\}\in\pi_1({\rm Ham}(M)),
  \end{equation}
 when $M$ is a Hirzebruch surface and when $M={\mathbb C}P^2$ (see Remark below Theorem \ref{Thmmain},
  the paragraph after
  Theorem \ref{Thmspheres} and Remark below Theorem \ref{ThmCPn}).


\medskip

{\bf Statement of main results}

\smallskip

{\bf 1. Necessary conditions for the $\sim$-equivalence.}
 An action
$\psi$ of $G$ on $M$ determines the fibration
$p:M_{\psi}:=EG\times_G M\to BG$, with fibre $M$, on the
classifying space $BG$ of the group $G$. If $\psi$ is Hamiltonian,
the fibration $M_{\psi}$ supports the coupling class $c_{\psi}\in
H^2(M_{\psi},\,{\Bbb R})$
 \cite{Mc-S}, \cite{G-L-S}. Then
fibre integration of   ${\rm exp}(c_{\psi})$ defines an element
$S_{\psi}$ in the cohomology $H(BG)$. We will prove that
$S_{\psi}=S_{\psi'},$ if the $G$-actions $\psi$ and $\psi'$ are
$\sim$-equivalent. That is,
\begin{Thm}\label{ThmHHA1} Let $\psi$ and $\psi'$ be two
Hamiltonian $G$-actions on $M$, if $S_{\psi}\ne S_{\psi'}$, then
$[\psi]\ne[\psi']\in[G,\,{\rm Ham}(M)]_{gh}$.
\end{Thm}

The cohomology class $S_{\psi}$  can be calculated by the
localization formula of the $G$-equivariant cohomology. This
formula gives an expression for $S_{\psi}$ which is a linear
combination of exponentials functions. If $G$ is the
$m$-dimensional torus, the coefficients of this combination are
rational functions of $m$ variables $(u_1,\dots, u_m)$.
 On the other hand, given  real constants
  $\gamma_{ij}$, the exponential functions
$F_i(u_1,\dots,u_m):={\rm exp}\big(\sum_j\gamma_{ij}u_j\big)$
  are ${\mathbb
R}(u_1,\dots,u_m)$-linearly independent (see Lemma
\ref{Lemexponbisbis}). This property will allow us to compare the
values of the map $S_{-}$ on different $G$-actions.
 The comparison of the values taken by this map, together with Theorem
 \ref{ThmHHA1}, will permit us to determine necessary
conditions,  with a simple geometric content,  for two torus
actions to be $\sim$-equivalent (see, for example, Theorem
\ref{ThmGEnralG} and Theorem \ref{Thmtorus}).

\smallskip

Let us assume that $(M,\omega)$ is equipped with a Hamiltonian
action of the Lie group $G$. If $X\in{\frak g}$ defines a circle
action $\psi^{X}$  on $M$, we say that $X$ {\em is of type $s$}
when a connected component of the fixed point set of $\psi^{X}$ is
a $2s$-dimensional submanifold of $M$, and there is no submanifold
of greater dimension in the fixed point set. The set of elements
of type $s$ will be denoted ${\rm Type}(s)$. We will prove the
following result, which gives a necessary condition for the
$\sim$-equivalence whose geometric meaning is clear

 \begin{Thm}\label{ThmGEnralG} If $X,Y\in{\frak g}$ are of different  type,
  then the corresponding circle actions $\psi^X$
 and $\psi^Y$ are not homotopic by means of a family of
 $U(1)$-actions.
 \end{Thm}

\smallskip

{\bf 2. Equivalence relations between torus actions.}

 In Section \ref{SectCohom} we will consider {\em toric
actions} on $M$; that is, Hamiltonian effective actions of the
$n$-dimensional torus $T:=U(1)^n$ on $M$. A toric action $\psi$
defines a  moment polytope $\Delta(\psi)=\Delta$ in ${\frak t}^*$,
and this polytope is uniquely determined by $\psi$ up to
translation.

Let $\psi,$ $\psi'$ be   toric actions on $M$, and let $\Delta,$
$\Delta'$ denote the moment polytopes  with center of mass at $0$
associated to $\psi$ and $\psi'$, respectively. By $\mu$ and
$\mu'$ we denote the corresponding moment maps, such that ${\rm
im}\,\mu=\Delta$ and ${\rm im}\,\mu'=\Delta'$.
 We state the
following properties:

 (i) $\,\Delta=\Delta'$.

(ii) $\,\Delta=\Delta'$ and there exists $h\in{\rm Symp}_0(M)$
such that $h\circ\psi_g\circ h^{-1}=\psi'_g$ for all $g\in T$,
where ${\rm Symp}_0(M)$ is the connected component of the identity
in the group ${\rm Symp}(M)$ of symplectomorphisms of $M$.





(iii) $\,\psi$ and $\psi'$ are  $r$-equivalent.

(iv) $\,S_{\psi}=S_{\psi'}$.

(v) $\,[\psi]=[\psi']\in[T,\,{\rm Ham} (M)]_{gh}.$

We will prove the following Theorem

\begin{Thm}\label{FinalThm}
The above properties satisfy the implications described in the
following diagram


$$\xymatrix{
 (iv)\ar@<0.5ex>[r] & (i)\ar@<0.5ex>[l] &\\
 (v)\ar[u]& (iii)\ar[l]\ar@<0.5ex>[r] &(ii)\ar@<0.5ex>[ul]\ar@<0.5ex>[l]
  \;.}
 $$

\end{Thm}

 From this theorem it follows,  among other properties,  that the
maps $S_{-}$ and $\Delta(\, .\,)$ are equivalent tools for
distinguishing $\sim$-inequivalent toric actions. The equality
$\Delta=\Delta'$ is a very simple geometric necessary condition
for the equality $[\psi]=[\psi']$,
  and a sufficient condition   can be obtained by
complementing that necessary condition as it is stated in (ii),
  but this complementary condition is not easy to check, in
  general.

\smallskip

 Given $\varphi,\,\varphi'\in{\rm Hom}(G,\,{\rm Ham}(M))$, we
  write $\varphi\sim_{rp}\varphi'$ if $\varphi$ and $\varphi'$ are $\sim$-equivalent ``up to
reparametrization" of $G$; that is, if there is an automorphism
$v$ of $G$ such that $\varphi\circ v\sim\varphi'$. Following
Karshon \cite{Kar} we say that the Hamiltonian actions $\varphi$
and $\varphi'$ are
 {\em conjugate} if they are  representations in  ${\rm Symp}(M)$ equivalent  ``up to
reparametrization"; that is, if there is $v\in{\rm Aut}(G)$ and
$\xi\in{\rm Symp}(M)$ such that $ \xi\circ(\varphi\circ
v)(g)\circ\xi^{-1}=\varphi'(g)$, for all $g\in G$. To complete the
scheme we define the $r$-equivalence up to reparametrization:
$\varphi$ and $\varphi'$ are $r$-equivalent up to
reparametrization iff there exist $v\in{\rm Aut}(G)$ and $h\in{\rm
Ham}(M)$ such that $h\circ(\varphi(v(g))\circ h^{-1}=\varphi'(g),$
for all $g\in G$.

Some of the implications in Theorem \ref{FinalThm} remains valid
when the equivalence relations are replaced by the corresponding
``up to  reparametrization" relations.
 For $\psi$ and $\psi'$ toric actions on $M$ we state the following
 properties:

 (a) $\psi$ and $\psi'$ are $r$-equivalent up to
 reparametrization.

 (b) $\psi\sim_{rp}\psi'.$

 (c) $\psi$ and $\psi'$ are conjugate.

 (d) There exists $V\in GL(n,\,{\mathbb Z})$, such that
 $V(\Delta)=\Delta'$.

  The following Proposition is also proved in Section 2

\begin{Prop}\label{Proeqiv}
For the above properties the implications of the following diagram

$$ \xymatrix{(a)\ar[r] \ar[d] &  (b) \ar[d] \\
 (c) \ar@<0.5ex>[r] & (d)\ar@<0.5ex>[l]    }$$
hold.

\end{Prop}

\smallskip

Let $M$ be  a toric manifold and $\psi$ the effective action of
$T$ which endows $M$ with the toric structure.
If $X$ belongs to the integer lattice of ${\frak t}$, it generates
an $S^1$-action $\psi^X$. If the element  $X$  is of type $s$,
then there is an
 $s$-face $f$ in the moment polytope $\Delta$ such that $\langle
 f,\,X\rangle=0$, but there is no an $r$-face orthogonal to $X$
 with $r>s$. If $X$ is of Type(0), the fixed point set for $\psi^X$ is the
 set $\{\mu^{-1}(P_i)\; |\; i=1,\dots, s\}$, where the $P_i$'s are the vertices of the moment polytope
 $\Delta$ and $\mu$ the moment map.

  We will prove the following Theorem that gives a simple
  geometric condition,
   under which two elements $X,Y$ of type $0$    define circle actions $\psi^X$ and $\psi^Y$
 no homotopic by a family of $U(1)$-actions.

 \begin{Thm}\label{Thmtorus}
 Let $M$ be a toric manifold and $\Delta$ the moment polytope with
  center of mass at $0$.
If $X,Y$ belong to the integer lattice of ${\frak t}$, they are in
${\rm Type}(0)$  and
$$ \{ \langle P_i,\,X\rangle\,:\,i=1,\dots,s \}\ne\{  \langle P_i,\,Y\rangle\,:\,i=1,\dots,s   \},  $$
$P_1,\dots,P_s$ being the vertices of $\Delta$. Then
$[\psi^X]\ne[\psi^Y]\in[U(1),\,{\rm Ham}(M)]_{gh}.$
 \end{Thm}

  Proposition \ref{Propalpha} in Section \ref{SectCohom}  is
the version ``up to reparametrization" of Theorem \ref{Thmtorus}.

\medskip

 {\bf 3. $\sim$-equivalence of $U(1)$-actions.}

\smallskip
  {\it $S^1$-actions on ${\mathbb P}(L_1\oplus\dots\oplus
L_{n-1}\oplus {\mathbb C})$.}

\smallskip

If $\Delta$ is a Delzant polytope in $({\mathbb R}^n)^*$, we
denote by $M_{\Delta}$ the symplectic toric  manifold associated
to $\Delta$ \cite{Del}, \cite{Gui}. We will study the relation
$\sim$ when $M_{\Delta}$ is diffeomorphic to the total space of
the fibration (\ref{DefPL}). In this case the result stated in
Theorem \ref{Thmtorus} will be completed giving a sufficient
condition for two elements in the integer lattice to define
$\sim$-equivalent $U(1)$-actions.

Let us assume that $M_{\Delta}$  is diffeomorphic to
 ${\mathbb P}(L_1\oplus\dots\oplus L_{n-1}\oplus{\mathbb
C})$, where $L_i$ is the holomorphic line bundle over ${\mathbb
 C}P^1$ with first Chern number $-a_i$. The $S^1$-action defined
 by ${\bf b}\in  {\mathbb Z}^n$ is denoted by $\psi_{\bf b}$.
 To study the relation $\sim$ between the $\psi_{\bf b}$'s we
 introduce the following definition:

 Given  ${\bf b}=(b_1,\dots,b_n)$ and ${\bf b'}=(b'_1,\dots,b'_n)$
two elements of the   lattice  ${\mathbb Z}^n$.
  We write ${\bf b}\equiv{\bf b'}$ iff either ${\bf b}={\bf
b'}$, or
$$ b'_j=b_j+a_jb_n,\;\;\text {for}\;\,j=1,\dots,
n-1 \;\;{\rm and}\;\,b'_n=-b_n.$$
 We will prove the
 following Theorem, which gives the aforesaid
 sufficient condition for the $\sim$-equivalence between $\psi_{\bf
 b}$ and $\psi_{\bf b'}$.
 \begin{Thm}\label{TmequivaPL}  Let ${\bf b}$ and ${\bf b'}$  be two elements of the integer
lattice   ${\mathbb Z}^n$. Let
  $\psi_{\bf b},$
$\psi_{\bf b'}$ be
  the respective circle actions on the symplectic   manifold
  $M_{\Delta}$ diffeomorphic to
  ${\mathbb P}(L_1\oplus\dots\oplus L_{n-1}\oplus{\mathbb C }).$
  If ${\bf b}\equiv{\bf b'}$,
then $\psi_{\bf b}\sim \psi_{\bf b'}$.
\end{Thm}

  When the manifold $M_{\Delta}$
is a Hirzebruch surface we will prove a stronger result than the
one stated in Theorem \ref{TmequivaPL}. Let $M_{\Delta}$ be the
Hirzebruch surface ${\mathbb P}(L_k\oplus {\mathbb C})$ (where
$L_k$ is the holomorphic line bundle over ${\mathbb C}P^1$ with
first Chern number $k$).
If ${\bf b}=(b_1,b_2)$ is an element in ${\mathbb Z}^2$, set
$S_{\bf b}=S_{\psi_{\bf b}}$.
 We will prove the following
Theorem
 \begin{Thm}\label{Thmmain}
 Let $M_{\Delta}$ be the Hirzebruch surface ${\mathbb P}(L_k\oplus {\mathbb
C})$, with $k\ne 0$, and ${\bf b},\, {\bf b'}\in {\mathbb Z}^2$,
then the following statements are equivalent

 (a) ${\bf b}\equiv {\bf b'}.$

(b) $\psi_{\bf b}\sim\psi_{\bf b'}.$

(c) $S_{\bf b}=S_{\bf b'}.$

\end{Thm}

\smallskip

{\it Remark.}
 Set ${\bf c}=(1,0)$ and  ${\bf c'}=(0,1)$. When $k$ is {\em odd},
  $$\{\psi_{\bf c}\},\,\{\psi_{\bf c'}\}\in \pi_1({\rm
Ham}({\mathbb P}(L_k\oplus {\mathbb C})))$$  are elements of
infinite order   in this group (see Section 4 of \cite{V1}). On
the other hand,  $\pi_1({\rm Ham}({\mathbb P}(L_k\oplus {\mathbb
C})))={\mathbb Z}$ (see \cite{Ab-McD}). Hence there exist
$m,n\in{\mathbb Z}\setminus\{0\}$ such that $\{\psi_{m{\bf
c}}\}=\{\psi_{n{\bf c'}}\}$. As the condition (a) of Theorem
\ref{Thmmain} does not hold for $m{\bf c}$ and $n{\bf c'}$, it
follows that $\psi_{m{\bf c}}$ and $\psi_{n{\bf c'}}$ are not
$\sim$-equivalent;  that is, the map (\ref{injective}) is not
injective for the Hirzebruch surface we are considering.

\smallskip

When $k=0$ the manifold ${\mathbb P}(L_k\oplus {\mathbb C})$ is
diffeomorphic to $S^2\times S^2$. In this case let  $\tau$ and
$\sigma$  denote the total areas of respective the spheres in the
product $S^2\times S^2$. We will prove following result

\begin{Thm}\label{Thmspheres}
If $b_1b_2\ne 0\ne b'_1b'_2$ and $\sigma/\tau\notin{\mathbb Q}$,
then the following statements are equivalent

 (i) $\; |b'_1|=|b_1|$ and $|b'_2|=|b_2|$.

 (ii) $\;S_{\bf b}=S_{\bf b'}.$

(iii)  $\;\psi_{\bf b}\sim \psi_{\bf b'}.$

  \end{Thm}

\smallskip

When $\sigma>1$ and $\tau=1$,
 \begin{equation}\label{pi1ham}
\pi_1({\rm Ham}(S^2\times S^2))={\mathbb Z}/2{\mathbb Z}\oplus
{\mathbb Z}/2{\mathbb Z}\oplus {\mathbb Z}.
\end{equation}
  The first two summands of the
right hand of (\ref{pi1ham}) come from the action of $SO(3)\times
SO(3)$ on $S^2\times S^2$. The generator of the third summand may
be represented by the loop $\{R^t \}_t$, where
$$R^t:(u,\,v)\in S^2\times S^2\mapsto (u,\,R^t_u(v))\in S^2\times S^2,$$
$R^t_u$ being the rotation of the fiber $S^2$, with angle $2\pi t$
and  axis through   $u$ (see \cite{Ab-McD}, \cite{McD}). Therefore
 $\psi_{\bf c}$, with $c=(1,0)$, is a representative of the generator of the first summand in the right
 hand of (\ref{pi1ham}).
  Thus $\{\psi_{\bf c}\}=\{\psi_{m\bf
c}\}\in\pi_1({\rm Ham}(S^2\times S^2))$, for all $m\in{2\,\mathbb
Z}+1$. But, according to Theorem \ref{Thmspheres}, from
$$[\psi_{\bf c}]=[\psi_{m\bf c}]\in[U(1),\, {\rm Ham}(S^2\times
S^2)]_{gh}$$
 it follows $m=\pm 1$, when $\sigma\notin{\mathbb Q}$.

\smallskip

  Although the results stated in Theorem \ref{Thmmain} and in
  Theorem \ref{Thmspheres},  the equality $S_{\bf b}=S_{\bf b'}$ is
  not equivalent to $[\psi_{\bf b}]=[\psi_{\bf b'}]\in
  [U(1),\,{\rm Ham}(M)]_{gh}$, in general. As we said, $S_{\bf b}=S_{\bf b'}$ is
  a necessary condition for the equality of the classes
  $[\psi_{\bf b}]$ and $[\psi_{\bf b'}]$. But if $(M,\omega)$ possesses certain symmetries
  it can happen that $S_{\bf b}=S_{\bf b'}$ and $[\psi_{\bf b}]\ne[\psi_{\bf
  b'}]$ (see Remark in Subsection \ref{Sectspheres}).

\medskip

{\it Circle actions on coadjoint orbits.}

  If $G$ is a compact connected Lie group and $\eta\in{\frak g}^*$, the stabilizer of $\eta$ for the coadjoint action
   contains a maximal torus $H$. The integer elements of ${\frak h}$ define Hamiltonian circle actions
    on the coadjoint orbit ${\mathcal O}$ of $\eta$.  We will consider the equivalence relation
 $\sim$ in the set of those  $S^1$-actions on the  orbit ${\mathcal
 O}$. By $W$ we   denote the Weyl group of the pair $(G,H)$. We
 will prove the following Theorem
\begin{Thm}\label{ThmCoadj}
Let $X, X'$ be two elements of   the integer lattice of ${\frak
h}$. If there exists $w\in W$ such that $w(X)=X'$, then the circle
actions defined on ${\mathcal O}$ by $X$ and $X'$ are
$\sim$-equivalent.
\end{Thm}

A particular case is when $G=SU(n)$, $H$ its diagonal subgroup and
${\mathcal O}$ is the orbit of an element whose stabilizer
contains $H$. In this case we have the following Corollary
\begin{Cor}\label{CoadjointSUn}
If $X=(x_1,\dots,x_n), \,X=(x_1,\dots,x_n)\in\big(2\pi i{\mathbb
Z})^n$ and there exists a permutation $\tau$ of the set
$\{1,\dots,n\}$, such that  $x'_j=x_{\tau(j)}$ for all $j$, then
the circle actions defined by $X$ and $X'$ on the orbit ${\mathcal
O}$ of $SU(n)$ are $\sim$-equivalent.
\end{Cor}

When the coadjoint orbit is diffeomorphic to ${\mathbb C}P^{n-1}$
we will prove a stronger result.


  \begin{Thm} \label{ThmCPn} Let ${\mathcal O}$ be a coadjoint orbit of $SU(n)$
  diffeomorphic to ${\mathbb C}P^{n-1}$, and $X,X'$ diagonal elements
  of ${\frak su}(n)$ such that $X, X'\in(2\pi i{\mathbb
  Z})^{n}$. If $\psi$ and $\psi'$ denote the circle actions on ${\mathcal O}$ defined by
 $X$ and $X'$ respectively. Then the following are equivalent

  (a) $[\psi]=[\psi']\in [U(1),\,{\rm Ham}(\mathcal O)]_{gh}.$

 (b) $S_{\psi}=S_{\psi'}.$

 (c) There exists a permutation $\tau$ of the set
$\{1,\dots,n\}$, such that  $x'_j=x_{\tau(j)}$ for  $j=1,\dots,n$.
 \end{Thm}

\smallskip
{\it Remark.} By the preceding Theorem, the elements $X\in (2\pi
i{\mathbb Z})^3$ define an infinite family in $[U(1),\,{\rm
Ham}({\mathbb C}P^2)]_{gh}$. But in contrast $\sharp(\pi_1({\rm
Ham}({\mathbb C}P^2)))=3$, since ${\rm Ham}({\mathbb C}P^2)$ has
the homotopy type of $PU(3)$ \cite{Gr}. Thus the map
(\ref{injective}) is not injective when $M={\mathbb C}P^2$.

\smallskip

To state a necessary condition for the equality $[\psi]=[\psi']$
when the coadjoint orbit of $SU(n)$ is the Grassmann manifold
$G_{k}({\mathbb C}^{n})$ (of the $k$-subspaces in ${\mathbb
 C}^n$) we need to introduce new notations.
 By ${\mathcal C}_k$ we denote the set of combinations of the set
 $\{1,\dots,n\}$ taken $k$ at a time. If $X=(ia_1,\dots,ia_n)$ and
${\frak d}=\{j_1,\dots,j_k  \}\in{\mathcal C}_k$ we write
$\alpha_{\frak{d}}=\sum_{l=1}^ka_{j_l}.$ We will prove the
following result

 \begin{Thm} \label{ThmGrass} Let ${\mathcal O}$ be a coadjoint orbit of $SU(n)$
  diffeomorphic to the Grassmann manifold $G_{k}({\mathbb C}^{n})$, and $X,X'$  diagonal regular elements
  of ${\frak su}(n)$ such that $X, X'\in(2\pi i{\mathbb
  Z})^{n}$. If $\psi$ and $\psi'$ denote the circle actions on $G_{k}({\mathbb C}^{n})$ defined by
 $X$ and $X'$ respectively. Then the equality
$$[\psi]=[\psi']\in [U(1),\,{\rm Ham}(\mathcal O)]_{gh}$$
implies that there exist a constant $\beta$ and  a bijective map
$f:{\mathcal C}_k\to {\mathcal C}_k$, such that $\alpha'_{\frak
d}=\alpha_{ f({\frak d})}+\beta$, for all ${\frak d}\in {\mathcal
C}_k$.

 \end{Thm}

{\it Remarks.}  When $k=1$ the
 necessary condition stated in
Theorem \ref{ThmGrass} for the equality $[\psi]=[\psi']$ is
 precisely the condition $(c)$ of Theorem \ref{ThmCPn}, since
 in this case the constant $\beta=a'_j-a_{f(j)}$ vanishes.

 On the other hand the mentioned necessary condition given in Theorem \ref{ThmGrass} is
independent of the symplectic structure. In Proposition
\ref{Fullflag} we give such a necessary condition  when the orbit
is diffeomorphic to the full flag manifold for $SU(n)$, which
depends on the symplectic structure.

\medskip

The paper is organized as follows. In Section \ref{SectCohom} we
introduce the cohomology class $S_{\psi}$ and we prove Theorem
\ref{FinalThm} and Proposition \ref{Proeqiv}.

Section \ref{SecPL} deals with the case when the manifold is the
total space of the fibration (\ref{DefPL}); in particular Theorem
\ref{TmequivaPL} is proved. Subsection \ref{Section4} is devoted
to prove Theorem \ref{Thmmain}.
 The case $M=S^2\times S^2$ is considered
in Subsection \ref{Sectspheres}.

In Section \ref{Coadjoint} we will concern with the relation
$\sim$ between circle actions  in the coadjoint orbits of a Lie
group $G$, defined by elements of a Cartan subalgebra of $G$.
We will prove the general Theorem \ref{ThmCoadj},  and Theorem
\ref{ThmCPn} and Theorem \ref{ThmGrass} relative to circle actions
on coadjoint orbits of $SU(n)$.

\smallskip

{\it Acknowledgements.} The  relation $\sim$   has been introduced
in \cite{V2}.
I thank to a referee of   paper \cite{V2} for suggesting me to
study the relation between the space $[G,\,{\rm Ham}(M)]_{gh}$ and
the $r$-equivalence of $G$-actions introduced above.  I thank the
referee  for his useful comments.


\section{The cohomology class $S_{\psi}$}\label{SectCohom}

  Here we will denote by ${\mathcal H}$   the Hamiltonian group ${\rm Ham}(M,\omega)$. From
   the universal principal bundle $E{\mathcal H}\to B{\mathcal H}$
   of the group ${\mathcal H}$, we construct the universal bundle
   with fibre $M$

$$
 \xymatrix{M\ar[r] & M_{\mathcal H}:=E{\mathcal H}\times_{\mathcal H}M
 \ar[d]^{\pi_{\mathcal H}} \\
 {} & B{\mathcal H}.}$$
There exists a unique class ${\bf c}\in H^2(M_{\mathcal H},\,{\Bbb
R})$  called the coupling class (see \cite{J-K}) satisfying:

a) $\pi_{{\mathcal H}*}({\bf c}^{n+1})=0$, where $\pi_{{\mathcal
H}*}$ is the fibre integration.

b) ${\bf c}$ extends the fiberwise class $[\omega]$.

Let $G$ be  a compact Lie group. We assume that $M$ is equipped
with the Hamiltonian $G$-action given by the group homomorphism
$\psi:G\to \text {Ham}\,(M,\,\omega)$. This homomorphism  induces
a map $\Psi:BG\to B{\mathcal H}$, between the corresponding
classifying spaces. The pullback $\Psi^{-1}(M_{\mathcal H})$ of
$M_{\mathcal H}$ by $\Psi$ is a bundle over $BG$ which can be
identified with $p:M_{\psi}=EG\times_G M\to BG$. Thus we have the
following commutative diagram
$$ \xymatrix{M_{\psi}\ar[r]^{\Hat\Psi }\ar[d]_p & M_{\mathcal H} \ar[d]^{\pi_{\mathcal H}}  \\
BG\ar[r]_{\Psi} & B{\mathcal H} \;.  }$$

By $c_{\psi}$ we denote the pullback of ${\bf c}$ by $\Hat\Psi$;
that is, $c_{\psi}=\hat\Psi^*({\bf c})$. The class $c_{\psi}$ is
the coupling class of the fibration $p: M_{\psi}\to BG$.

 If $\psi$ and
$\psi'$ are homotopic by means of a family of Hamiltonian
$G$-actions, then the bundles $\Psi'^{-1}(M_{\mathcal H})$ and
$\Psi^{-1}(M_{\mathcal H})$ are isomorphic and the isomorphism
from $M_{\psi}$ to $M_{\psi'}$ applies $c_{\psi'}$ to $c_{\psi}$.

Since $c_{\psi}\in H^2(M_{\psi})=H^2_G(M)$ we can integrate
$e^{c_{\psi}}$ along the fibre and we obtain
$S_{\psi}:=p_*(e^{c_{\psi}})$ in the cohomology of $H(BG)$. {\em
The above arguments prove   Theorem \ref{ThmHHA1}.}

\smallskip

Strictly speaking $e^{c_{\psi}}$ is not an equivariant de Rham
cohomology class. For instance, when $G$ is the $n$-torus $T$,
$e^{c_{\psi}}$ belongs to $\Hat H_T(M)$ \cite{A-B}, where $\Hat
H_T(M)$ is
$$\Hat H_T(M)=H_T(M)\otimes_{{\Bbb R}[u_1,\dots,u_n]}{\Bbb R}[[u_1,\dots,u_n
]],$$
  ${\Bbb R}[[u_1,\dots,u_n ]]$ being  the ring of formal power
  series in the variables $u_1,\dots, u_n$.

 By the localization theorem \cite{A-B} we have, in a suitable
localization of $\Hat H(BT)={\Bbb R}[[u_1,\dots, u_n]]$,
\begin{equation}\label{Localiz}
 S_{\psi}=\sum_Fp_*^F\Big(\frac{e^{c_{\psi}}}{E(N_F)} \Big),
\end{equation}
 where $F$ varies in the set
of connected components of the fixed point set of the $T$-action,
$p^F_*:\Hat H_T(F)\to \Hat H(BT)$ is the fibre integration on $F$,
and $E(N_F)$ is the equivariant Euler class of the normal bundle
$N_F$ to $F$ in $M$.

\smallskip

 Let us assume that $M$ is endowed with a Hamiltonian action of the torus $T=U(1)^n$.
  We say that a moment map $\mu:M\to {\frak
t}^*$ for the $T$-action is the {\em normalized} moment map if
$\int_M \mu(X)\omega^n=0$, for all $X\in{\frak t}$. If $\Delta$
denotes the moment polytope of $\mu$, that is $\Delta=\rm
Im(\mu)$, then
$$\int_M\mu\,\omega^n=\text {Cm}(\Delta)\int_M\omega^n,$$
where $\rm Cm(\Delta)$ is the center of mass of $\Delta$. It
follows that $\rm Cm(\Delta)=0$, iff $\mu$ is normalized.

If $\mu$ is a normalized moment for the $T$-action $\psi$, then
$\omega-\mu$ is a $T$-equivariant closed $2$-form \cite{G-S} which
represents the coupling class $c_{\psi}$, because it satisfies the
properties which determine the coupling class uniquely.

Let $M$ be  a toric manifold and $\psi$ the corresponding
effective action of $T$. We denote by $\Delta$ the respective
moment polytope with center of mass at $0$. By $u_1,\dots, u_n$ we
denote the basis of ${\frak t}^*$ dual of $e_1,\dots, e_n\in
{\frak t}={\Bbb R}^n$. The image by moment map $\mu$ of the fixed
points $q_1,\dots, q_s$ of the $T$-action $\psi$ are the vertices
of $\Delta$ \cite{At}, \cite{vG82}. The equivariant Euler class
$E(N_{q_i})$ of the normal bundle to $q_i$ in $M$ is
$$E(N_{q_i})=\prod _{j=1}^nl_{ij},$$
where the $l_{ij}$ are the weights of isotropy   representation of
$T$ at the tangent space $T_{q_i}$ (see \cite{G-S}). So
$E(N_{q_i})$ is a homogeneous degree $n$ polynomial in the
variables $u_1,\dots, u_n$.

On the other hand $\mu(q_i)=\sum_j\beta_{ij}u_j$, where
$(\beta_{ij})_j$ are the coordinates of the vertex $\mu(q_i)$. Let
$f$ be the polynomial $\prod_i\prod_j l_{ij}$. Denoting by
$[S_{\psi}]_f$ the class of $S_{\psi}$   in the localization $\Hat
H(BT)_f$ (see \cite{A-B}) it follows from (\ref{Localiz})
\begin{equation}\label{Loocalim}
[S_{\psi}]_f=\sum_{i=1}^s\frac{e^{-\sum_j \beta_{ij} u_j}}{\prod_j
l_{ij}}\in \Hat H(BT)_f.
\end{equation}
So we can write
\begin{equation}\label{Spsif}
 [S_{\psi}]_f=\frac{A}{f}\in \Hat H(BT)_f,
\end{equation}
where  $A$ is  a linear combination of elements of
 $${\mathcal
E}:=\{E_i:={\rm  exp}(-\sum_j \beta_{ij}u_j)\}_i$$
 with
coefficients in the polynomial ring ${\Bbb R}[u_1,\dots, u_n]$.

 Let $\psi'$ be another effective action of $T$ on $M$, and  $f'$
 be
 the polynomial defined by  the equivariant Euler classes of the normal bundles of
 the fixed points. So
\begin{equation}\label{Spsif'}
 [S_{\psi'}]_{f'}=\frac{C}{f'}\in\Hat H(BT)_{f'}
\end{equation}

In order to compare $[S_{\psi}]_f$ and $[S_{\psi'}]_{f'}$ we
introduce the extension ${\mathcal A}$ of $\Hat H(BT)$
$${\mathcal A}:={\Bbb R}[[u_1,\dots,u_n]]\otimes_{{\Bbb
R}[u_1,\dots,u_n]}{\Bbb R}(u_1,\dots, u_n),$$
 where ${\Bbb R}(u_1,\dots, u_n),$ is the field of rational
 functions in the variables $u_1,\dots,u_n$.

We have the following obvious Lemma

 \begin{Lem}\label{Lema1} Suppose that $S_{\psi}=S_{\psi'}\in
 \Hat H(BT)$. If $[S_{\psi}]_{f}=\frac{\tilde A}{f^p}\in\Hat
 H(BT)_{f}$ and
 $[S_{\psi'}]_{f'}=\frac{\tilde C}{(f')^q}\in\Hat H(BT)_{f'}$, with $p,q\in{\mathbb Z}$,  then there exists a rational function $g\in {\Bbb
 R}(u_1,\dots,u_n)$ such that $\tilde A=\tilde{C} g$ in ${\mathcal A}$.
\end{Lem}
{\it Proof.}   $g=f^p/(f')^q$ satisfies the required condition.
\qed

\smallskip

Let $\{\gamma_{ij}\,|\, i=1,\dots, r; \,j=1,\dots, m  \}$ be a
family of real numbers. We put
$$F_i(u_1,\dots,u_m)={\rm exp}\big(\sum_j\gamma_{ij}u_j\big).$$
 We
will prove that the set of exponential functions $\{F_i\}\in
{\mathcal A}$ is ${\mathbb R}(u_1,\dots,u_m)$-independent.
 \begin{Lem}\label{Lemexponbisbis}
 If $\{Q_i(u_1,\dots,u_m)\}_{i=1,\dots,r}$ is a set of rational
 functions, such that $\sum_iQ_iF_i=0$ and $F_i\ne F_k$, for $i\ne
 k$, then
 $Q_i=0$ for $i=1,\dots, r.$
\end{Lem}

{\it Proof.} We define the linear functions
$h_i(u_1,\dots,u_m)=\sum_j\gamma_{ij}u_j.$ We write
$|\gamma_i|=(\sum_j\gamma_{ij}^2)^{1/2}.$ Let us assume that
$|\gamma_1|\geq |\gamma_i|$ for $i=1,\dots, r$. Then $|{\rm
grad}(h_i)|\leq |{\rm grad}(h_1)|$ and moreover ${\rm
grad}(h_i)\ne {\rm grad}(h_1)$, for $i\ne 1$ by the assumption.
Hence the vector $v:=(\gamma_{11},\dots,\gamma_{1m})\in{\mathbb
R}^m$ defines a direction such that
\begin{equation}\label{Quotienlimi}
\lim_{\lambda\to+\infty}\frac{F_i(\lambda v)}{F_1(\lambda v)}=0,
\end{equation}
for $i\ne 1$.

If $Q_1$ would be nonzero
\begin{equation}\label{Equ1=0}
1=-\sum_{i>1}\frac{Q_i}{Q_1}\frac{F_i}{F_1}.
\end{equation}
 If we
take the limit along $\lambda v$ as $\lambda\to+\infty$ in
(\ref{Equ1=0}) we arrive to a contradiction, by the exponential
decay (\ref{Quotienlimi}). Hence the initial linear combination
reduces to $\sum_{i>1}Q_iF_i=0$. By repeating the argument it
follows $Q_i=0$ for all $i$.
 \qed

\smallskip

Sometimes we will use particular cases of Lemma
\ref{Lemexponbisbis}; for example,  the ${\mathbb
R}[u,u^{-1}]$-independence of a family of exponential functions
$\{e^{\gamma_i u}\}_{i}$,  with $\gamma_i\ne\gamma_j$ for $i\ne
j$.

 \begin{Prop}\label{CenPro}
 Let $\psi$ and $\psi'$ be   toric actions on $M$ and $\Delta$ and $\Delta'$ be the respective
  moment polytopes with center of mass at $0$. Then
  $\Delta\ne\Delta'$ iff $S_{\psi}\ne S_{\psi'}$
  \end{Prop}
{\it Proof.}  Let us suppose that $\Delta\ne \Delta'$. According
to (\ref{Loocalim}),
 the element $A$ in (\ref{Spsif}) is a linear combination
 $A=\sum_ig_iE_i,$ with
  $0\ne g_i\in{\Bbb R}[u_1,\dots,u_n].$
  Similarly,
$C$ in (\ref{Spsif'}) is a linear combination of elements of
 $${\mathcal
E'}:=\{E'_i:={\rm  exp}(-\sum_j \beta'_{ij}u_j)\}_i.$$
  As the
moment polytopes are not equal, there is a vertex in $\Delta$
which does not belong to $\Delta'$. So there is an element, say
$E_k$,  in ${\mathcal E}$ which is not in ${\mathcal E'}$.

 If there were a rational function $g\in{\Bbb R}(u_1,\dots,u_n)$ such that $A=Cg$ in ${\mathcal A}$, then
$0=A-Cg$ is a linear combination in ${\mathcal A}$
 with coefficients in ${\Bbb R}(u_1,\dots,
u_n)$, in which  $E_k$ appears  once with the coefficient $g_k\ne
0.$  That is,
\begin{equation}\label{auxequa}
0=A-gC=g_kE_k+\sum_{G_r\in{\mathcal E"}} \tilde{g}_rG_r,
\end{equation}
where ${\mathcal E"}=({\mathcal E}\cup {\mathcal E'})\setminus
\{E_k\}$. By the ${\Bbb R}(u_1,\dots,u_n)$-linear independence of
the set ${\mathcal E"}\cup \{ E_k \}$, it follows from
(\ref{auxequa}) that $g_k=0$.
This contradiction shows that there
 there is no a rational function $g$ such that $A=Cg$. It follows from
 Lemma
\ref{Lema1} that $S_{\psi}\ne S_{\psi'}$.

On the other hand, if $\Delta=\Delta'$, then trivially
$\beta_{ij}=\beta'_{ij}$. Moreover the equality of the polytopes
implies the equality of the weights of the respective isotropy
representations at the fixed points \cite{vG82}. Hence, by
(\ref{Loocalim}), $[S_{\psi}]_f=[S_{\psi'}]_{f'}$. As the
polynomial $f$ and $f'$ are equal, it turns out that
$S_{\psi}=S_{\psi'}$.

  \qed

\begin{Thm}\label{Delz}
Let $\psi$ and $\psi'$ be  toric actions on $M$. Let $\Delta$ and
$\Delta'$ be the respective
  moment polytopes with center of mass at $0$. If
  $\Delta\ne\Delta'$, then
  $[\psi]\ne[\psi']\in[T,\,\rm
  Ham(M)]_{gh}$.
\end{Thm}

{\it Proof.} It is a consequence of Proposition \ref{CenPro} and
Theorem \ref{ThmHHA1}.
 \qed

\begin{Prop}\label{Propref} If $\psi$ and $\psi'$ are $r$-equivalent Hamiltonian
$G$-actions on $M$, then $[\psi]=[\psi']\in[G,\,{\rm
Ham}(M)]_{gh}$.
\end{Prop}

{\it Proof.} There exists $h \in{\rm Ham}(M)$ such that
$h\circ\psi_g\circ  h^{-1}=\psi'_g$, for all $g$.
 Let $h_s$ be a path in ${\rm Ham}(M)$ from
${\rm Id}$ to $h$, then $\{h_s\circ\psi\circ h_s^{-1}\}_s$ defines
a homotopy of $G$-actions between $\psi$ and $\psi'$.
 \qed

\medskip

{\bf Proof of Theorem \ref {FinalThm}.}

 Properties (iv) and  (i) are equivalent by Proposition
\ref{CenPro}.

 (v) ``implies" (iv)  by Theorem \ref{ThmHHA1}.

 (iii) $\,\Longrightarrow\,$ (v) is
 Proposition \ref{Propref}.

 If $h\in{\rm Ham}(M)$ satisfies $\psi'_g=h\circ\psi_g\circ h^{-1}$, for all $g\in T$, then the corresponding
 normalized moment maps satisfy $\mu'=\mu\circ h $. So
$\Delta={\rm im}\,\mu={\rm im}\,\mu'=\Delta'$; that
 is, (iii) implies (ii).

 Since $M$ is a toric manifold ${\rm Ham}(M)={\rm Symp}_0(M)$; thus (ii)
 implies (iii).
\qed

\begin{Lem}\label{Lemat} Let $\psi$ be a toric action on $M$ and $\Delta$ its moment polytope with center of mass at
$0$. If    $v$ is an automorphism of $T$, then there exists $V\in
GL(n,\,{\Bbb Z})$ such that $V(\Delta)$ is the moment  polytope
with center of mass at $0$ associated to the action $\psi\circ v$.
\end{Lem}
{\it Proof.} Let $v^*$ be the automorphism of ${\frak t}^*$
induced by $v$. If $\mu$ is the normalized moment map of $\psi$,
then $v^*\circ\mu$ is the corresponding moment map of $\psi\circ
v$. So $V=v^*$ in the identification ${\frak t}^*=({\mathbb
R}^n)^*$.

\qed

\smallskip

{\bf Proof of Proposition \ref{Proeqiv}.}

(a)``implies" (b). It is consequence of Proposition \ref{Propref}.

 (b)``implies" (d). By definition of $\sim_{rp}$,
there exists  an automorphism $v$ of $T$ such that $\psi'\sim
\psi\circ v$. By the implication (v)$\,\Longrightarrow\,$ (i) of
Theorem \ref{FinalThm}, $\Delta'=\Delta(\psi\circ v)$. The
implication follows from Lemma \ref{Lemat}.

(c) ``implies " (d). If $\psi'$ and $\psi$ are conjugate then
$\xi\circ(\psi\circ v)(g)\circ\xi^{-1}=\psi'(g)$, with $\xi$ a
symplectomorphism. By the proof of Lemma \ref{Lemat},
$\mu'=v^*\circ\mu\circ\xi$. So $\Delta'=v^*(\Delta)$.

(d) ``implies " (c). Let us assume that $V(\Delta)=\Delta'$ with
$V\in GL(n,{\mathbb
 Z})$. $\,V$ induces a group automorphism of $T$, which is denoted
 by $v$, such that $v^*=V$. So the actions $\psi\circ v$ and
 $\psi'$ have the same moment polytope. By
 Delzant's Theorem (Th. 2.1, \cite{Del}) there is a symplectomorphism $\xi$ such that
 $\xi\circ(\psi\circ v)(g)\circ\xi^{-1}=\psi'(g)$; that is, $\psi$
 and $\psi'$ are conjugate.

 \qed

\smallskip

Let us assume that $M$ is a symplectic   manifold   equipped with
a Hamiltonian $G$-action, where $G$ is a Lie group. Let $X$ be an
element of ${\frak g}$ of Type (s) which determine a circle action
$\psi^X$ on $M$, and let $F$ be a $2s$-dimensional submanifold of
$M$, connected component of the fixed point set of the circle
action $\psi^X$. The normal bundle $N_F$ to $F$ in $M$ decomposes
in direct sum of $n-s$ $U(1)$-equivariant line bundles
$$N_F=\bigoplus_{j=1}^{n-s} N_j.$$
  On $N_k$ $\,U(1)$ acts with a weight  $p_k u\ne 0$,
  where $u$ we denote the   vector in ${\frak u}(1)^*$ dual of the
standard basis of ${\frak u}(1)$.
  Then the $U(1)$-equivariant Euler class of $N_{F}$ is
  $$E(N_F)=\prod_{k=1}^{n-s}\big( c_1(N_k)+p_k u \big)
   =u^{n-s}\prod_{k} p_k\prod_{k }\Big(1+\frac{c_1(N_k)}{p_k u}
   \Big).$$
   So
\begin{equation}\label{Eqparadespues}
 \big(E(N_F)\big)^{-1}= u^{s-n}\Big(\alpha_0+\frac{\alpha_1}{u}
+\frac{\alpha_2}{u^2}+\dots \Big),
 \end{equation}
 $\alpha_r$ being a $2r$-form on $F$ and $\alpha_0\ne 0$. That is,
$(E(N_F))^{-1}$ is a Laurent polynomial of
 degree $s-n$ with coefficients even dimensional forms on $F$.

 On the other hand, if $f$ is a normalized Hamiltonian for the $U(1)$-action $\psi^X$, then  $\omega -uf$
 is a $U(1)$-equivariant closed $2$-form \cite{G-S} which
 represents the coupling class $c_{\psi^X}$  \cite{K-M}.
Hence the contribution of $F$ to $S_{\psi_X}$
 in (\ref{Localiz})  is
 $$ \frac{e^{au}}{u^{n-s}}\int_Fe^\omega\Big( \alpha_0+\frac{\alpha_1}{u}
 +\dots \Big),$$
  where $-a$ is the value of $f$ on $F$. That is,
$$ \frac{e^{au}}{u^{n-s}}  \Big( A_0 u+ A_1+\dots +\frac{A_{m-1}}{u^{m-2}}
  \Big),$$
 with $A_0\ne 0$.

 Let $F'$ be another component of the fixed point set.
 As ${\rm codim}\,F'\geq {\rm codim}\,F\geq 2$, there is a curve
 $\sigma$ in $M$ which joins $F$ and $F'$ and such that it does
 not meet  other fixed point. So the only stationary points of the
 function $f\circ \sigma$ are the end points; that is, $f\circ
 \sigma$ is strictly monotone, and thus $f(F)\ne f(F')$. Therefore
 the contribution of $F'$ to $S_{\psi^X}$ will be of the form
 $e^{a'u}Q'(u)$, with $a\ne a'$ and $Q'(u)$ a Laurent polynomial
 of degree $(1/2)\,{\rm dim}\,F'-n$. We note that the
 contributions of $F$ and $F'$
 can not be grouped together in the form $e^{\tilde a
 u}\tilde Q(u)$. So we have the following Lemma
 \begin{Lem}\label{Lemgene}
If $X$ is of Type (s) and $F_1,\dots,F_c$ are the $2s$-dimensional
components of the fixed point set of the circle action $\psi^X$,
then
$$S_{\psi^X}=\sum_{i=1}^c e^{a_i u}Q_i(u)+\sum_r e^{b_r u}\Hat
Q_r(u),$$ with $a_i\ne a_j$ for $i\ne j$, $Q_i(u)$ a Laurent
polynomial of degree $s-n$ and $ \Hat Q_r(u)$ a Laurent polynomial
of degree less than $s-n$.
 \end{Lem}

\smallskip

{\bf Proof of Theorem \ref{ThmGEnralG}.}  If $X$ and $Y$ are of
different Type, say Type of $X=s$ greater than Type of $Y$, then
in the expression for $S_{\psi^X}$ given in Lemma \ref{Lemgene}
appears an exponential function with coefficient a Laurent
polynomial of degree $s-n$, but all the Laurent polynomials in
$S_{\psi^Y}$ have degree less than $s-n$. The Theorem follows from
 the aforesaid ${\mathbb R}[u,u^{-1}]$-linear  independence of
the exponential functions $\{e^{au}\}$ together with Theorem
\ref{ThmHHA1}.
  \qed

\smallskip

  Let $M$ be a toric manifold. To avoid trivial cases we
assume that ${\rm dim}\,M>2$.  By $\psi$ we denote the
corresponding Hamiltonian effective $T$-action
and by $\mu$ a
 moment map. In terms of the basis $\{u_j\}$ we write
$\mu=\sum_j \mu_ju_j$. Given $X=\sum_jk_je_j$ in the integer
lattice of ${\frak t}$, then $\mu^X=\langle\mu,\,X\rangle=\sum_j
k_j\mu_j$ is a moment map for the $U(1)$-action $\psi^X$ defined
as the composition of $\psi$ with $z\in U(1)\mapsto (z^{k_j})\in
T$.

The following Proposition gives equivalent conditions to the
property of being   of type $0$.

\begin{Prop}\label{Proequ}
Let $X$ be a vector in the integer lattice of ${\frak t}$, the
following are equivalent.

($\alpha$) The set of fixed points for the  $U(1)$-action $\psi^X$
is equal to the set of fixed points of the $T$-action $\psi$.

($\beta$) $\,\langle v, X\rangle\ne 0,$ for all $v$ parallel to an
edge of the moment polytope $\Delta$.

($\gamma$) $\,X\notin{\frak g}_a={\rm Lie}(G_a),$  for all $a$,
with $G_a$ proper isotropy subgroup of the $T$-action $\psi$.

\end{Prop}

{\it Proof.}

 ($\alpha$)  $\,\Longrightarrow\,$  ($\beta$). Let $e$ be an edge  of moment
polytope $\Delta$. $e$ is of the form $p+tv$, with $t$ in an
interval of ${\Bbb R}$, $p$ a vertex of $e$ and $v\in{\frak t}^*$.
If $\langle v, X\rangle=0$, then $\mu^X$
  takes on $\mu^{-1}(e)$ the value $\langle
p+tv,\,X\rangle= \langle p,X\rangle$. Thus the infinitely many
points of $\mu^{-1}(e)$ are fixed by $\psi^X$. But the fixed point
set of the $T$-action is $\{q_1,\dots,q_s\}$, the set of inverse
images by $\mu$ of the vertices of $\Delta$.

 ($\beta$)   $\,\Longrightarrow\,$  ($\gamma$). If $X\in{\frak g}_b$, we may
assume that ${\frak g}_b$ is maximal; that is, $G_b$ is the
stabilizer group of the points of $N:=\mu^{-1}(e)$, with $e$ an
edge of $\Delta$. So the induced vector field $X_M$ on $M$
vanishes on the submanifold $N$, and $d\mu^X=0$ on this set. The
values of the function $\mu^X$ on $N$ are of the form $\langle
p+tv,\,X\rangle$, with $t$ varying in an interval of ${\Bbb R}$,
and where $p$ is a vertex of $e$ and $v$ is parallel to $e$. It
follows from the vanishing of $d\mu^X|_{N}$ that $\langle v,
X\rangle=0$.

 ($\gamma$)  $\,\Longrightarrow\,$ ($\alpha$). If $X$ satisfies ($\gamma$) then the
fixed point set for the $U(1)$-action $\psi^X$ is
$\{q_1,\dots,q_s\}$ (see Corollary 10.9.1 in \cite{G-S}).

\qed

If the equivalent properties stated in Proposition \ref{Proequ}
hold, the fixed point set for $\psi^X$ is $\{q_1,\dots, q_s \}$.
The $U(1)$-equivariant
 Euler class
 $E(N_{q_i})$ is a degree $n$ monomial in a variable $u$, that we
 denote by  $l_i$. If the moment map $\mu$ is normalized, then
 $\mu^X$ is the normalized Hamiltonian of the $S^1$-action
 $\psi^X$.
 Therefore
 $$[S_{\psi^X}]=\sum_{i=1}^s\frac{e^{-\beta_i u}}{l_i}\in \Hat
 H(BU(1))_u,$$
 with $\beta_i=\langle P_i,\,X\rangle,$ and $P_i$ the vertex of
 $\Delta$ image of $q_i$ by $\mu$.
 The proof of    Theorem \ref{Thmtorus} is similar to the one given
 for Proposition \ref{CenPro}.

\smallskip

 {\bf Proof of Theorem \ref{Thmtorus}.} We put
${\mathcal E}=\{E_i=e^{-\beta_i u} \}$ and ${\mathcal
E'}=\{E'_i=e^{-\beta'_i u}  \}$, with $\beta'_i=\langle
\mu(q_i),\,Y\rangle$. We write $[S_{\psi^X}]=A/u^n$ and
$[S_{\psi^Y}]=C/u^n$, where $A,C\in{\Bbb R}[[u]]$ are linear
combinations
$$A=\sum_{i=1}^sa_iE_i,\;\;\;C=\sum_{i=1}^sc_iE'_i;\;\;\;
a_i\ne 0\ne c_i, \;\;\; {\rm for all}\;\; i.$$
 By hypothesis,
there is an $E_k$ that does not belong to ${\mathcal E'}$. Then,
by the ${\Bbb R}$-linear independence of the exponential
functions, $A\ne  C$.
  Thus $S_{\psi^X}\ne S_{\psi^Y}$, and by Theorem
\ref{ThmHHA1} $[\psi^X]\ne[\psi^Y]\in [U(1),\,{\rm Ham}(M)]_{gh}.$

\qed

\begin{Prop}\label{Propalpha}
Let $X,Y$ be elements in the integer lattice of ${\frak t}$ that
satisfy the conditions of Proposition \ref{Proequ}. If
 \begin{equation}\label{auxPropalpha}
 \{\lambda
\langle P_i,\,X\rangle\,:\,i=1,\dots,s \}\ne\{  \langle
P_i,\,Y\rangle\,:\,i=1,\dots,s   \},
 \end{equation}
   for all $0\ne\lambda\in{\mathbb Z}$, where the $P_i$'s are the
   vertices of the moment polytope with center at $0$.
Then $\psi^X$ and
 $\psi^Y$ are $S^1$-actions $\sim_{rp}$-inequivalent.
\end{Prop}

{\it Proof.} Suppose the Proposition were false. Then there were
an automorphism $v$ of $U(1)$ such that $\psi^X\circ v\sim\psi^Y$.
But $\psi^X\circ v=\psi^Z$, with $Z=v_*(X)$; that is, $Z=\lambda
X$, with $0\ne \lambda\in{\mathbb Z}$. Since $\langle
P_i,\,Z\rangle=\lambda\langle P_i,\,X\rangle$, it follows from
(\ref{auxPropalpha}) together with Theorem \ref{Thmtorus} that
$[\psi^Z]\ne[\psi^Y]$; that is, $\psi^X\circ v$ is not
$\sim$-equivalent to $\psi^Y$. This  contradicts  our supposition.

\qed

\section{Circle actions  in  ${\mathbb P}(L_1\oplus\dots\oplus L_{n-1}\oplus{\mathbb C
})$.}\label{SecPL}

We denote by  $e_1,\dots,e_n$ the standard basis of ${\mathbb
R}^n$.  Given $a_1,\dots,a_{n-1}\in{\mathbb Z}$,
 we define the following vectors of ${\mathbb R}^n$:
 $$v_j=-e_j,\;
j=1,\dots,n, \;\;\; v_{n+1}=\sum_{j=1}^{n-1}e_j,\;\;\;
v_{n+2}=e_n-\sum_{i=1}^{n-1}a_ie_i.$$
 Let $\{k_i\}_{i=1,\dots, n+2}$ be a family of real numbers
with $k_j=0$, for $j=1,\dots, n$, $k_{n+1}=\sigma>0$,
$k_{n+2}=\tau>0$ satisfying
  $\tau+a_i\sigma>0$ for all $i$, then  we consider the following
  Delzant polytope $\Delta$ in  $({\mathbb R}^n)^*$
$$\Delta=\bigcap_{i=1}^{n+2}\{x\in ({\mathbb R}^n)^*\,:\, \langle x,\,v_i\rangle\leq k_i  \}.$$
  When $n=3$ the polytope $\Delta$ is represented in
the Figure 1.

\begin{figure}[htbp]
\begin{center}
\epsfig{file=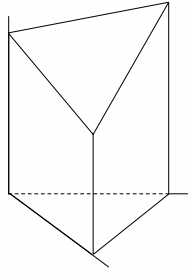, height=4cm}
\end{center}
\caption[Figure 1]{\small The polytope $\Delta$ when $n=3$}
\end{figure}

The symplectic manifold $M_{\Delta}$ associated to $\Delta$ is
$$M_{\Delta}=\{ (z_i)\in {\mathbb C}^{n+2}\,:\, \sum_j\, '\,|z_j|^2=\sigma/\pi,\;
-\sum_{j=1}^{n-1}a_j|z_j|^2+|z_n|^2+|z_{n+2}|^2=\tau/\pi
\}/\simeq,$$ in the sum $\sum_j '$ the index $j$ runs in the set $
\{1,\dots,n-1,n+1\}$. The equivalence relation $\simeq$ is defined
by $(z'_j)\simeq (z_j)$ iff there exist $\epsilon,\lambda\in U(1)$
such that
$$\begin{cases}z'_j=\lambda\epsilon^{-a_j}z_j,\; j=1,\dots,n-1;\\
z'_k=\epsilon z_k,\; k=n,\,n+2;\\  z'_{n+1}=\lambda z_{n+1}.
\end{cases}$$

The action of $U(1)^n$ on $M_{\Delta}$ given by
\begin{equation}\label{TAction}
 (\lambda_j)\cdot[z_i]=[\lambda_1z_1,\dots,\lambda_nz_n,z_{n+1},z_{n+2}]
\end{equation}
defines a structure of toric manifold on $M_{\Delta}$. Moreover
the map $\mu$ on $M_{\Delta}$ defined by $\mu([z_j])=
\pi(|z_1|^2,\dots,|z_n|^2)$ is the corresponding moment map with
${\rm im}\,\mu=\Delta.$

The manifold $M_{\Delta}$ is the total space of the fibration
${\mathbb P}(L_1\oplus\dots\oplus L_{n-1}\oplus{\mathbb C })\to
{\mathbb C}P^1$, where $L_i$ is the holomorphic line bundle over
${\mathbb C}P^1$ with first Chern number $-a_i$.

Given $[z]\in M_{\Delta}$, we write $w=(z_n,\,z_{n+2})$. The
element in ${\mathbb P}(L_1\oplus\dots\oplus L_{n-1}\oplus{\mathbb
C })$ which corresponds to $[z]$ is the following one in the fibre
over $[w]=[z_n:\,z_{n+2}]\in{\mathbb C}P^1$
\begin{equation}\label{ElemenPL}
[(z_n)^{a_1}z_1:\,(z_{n+2})^{a_1}z_1:\,\dots :\,
(z_n)^{a_{n-1}}z_{n-1}:\,(z_{n+2})^{a_{n-1}}z_{n-1}: \,z_{n+1}  ]
\end{equation}
We will write $([w],\,[w^{a_j}z_j:\,z_{n+1}])$ for the element in
fibre bundle associated to $[z]\in M_{\Delta}$.

An element ${\bf b}=(b_1,\dots,b_n)$ in the integer lattice
${\mathbb Z}^n$   determines a Hamiltonian circle action
$\psi_{\bf b}$ on $M_{\Delta}$ through the  action
(\ref{TAction}). In the notation introduced above
$$\psi_{\bf b}(t)[z]=\big( [w_t],\, [(w_t)^{a_j}z_j\epsilon_j:\,z_{n+1}]
\big),$$
 where $\epsilon_k={\rm exp}(2\pi ib_k t)$ and
$w_t=(z_{n}\epsilon_{n},\,z_{n+2}).$

\smallskip

{\bf Proof of Theorem \ref{TmequivaPL}.} Let ${\bf
b'}=(b'_1,\dots,b'_n)$ with
$$ b'_j=b_j+a_jb_n,\;\;\text {for}\;\,j=1,\dots,
n-1,\;\;{\rm and}\;\,b'_n=-b_n,$$
 then
$$\psi_{\bf b'}(t)[z]=\big( [w'_t],\, [(w'_t)^{a_j}z_j\epsilon'_j:\,z_{n+1}]
\big),$$
 with $w'_t=(z_{n}\epsilon'_{n},\,z_{n+2})$ and $\epsilon'_k={\rm exp}(2\pi ib'_k t)$. So
$$\psi_{\bf b'}(t)[z]=\big( [u_t],\, [(u_t)^{a_j}z_j\epsilon_j:\,z_{n+1}]
\big),$$
  where $u_t=(z_{n},\,z_{n+2}\epsilon_{n}).$

 Let $\{\Hat n_s\,:\, s\in[0,1]\}$ be a continuous curve of unitary
 vectors in ${\mathbb R}^3$, with $\Hat n_0= e_3$ and $\Hat n_1=- e_3$. We denote by
   $\xi^s(t)$  the rotation of $S^2={\mathbb C}P^1$
around the axis defined by the unitary vector $\Hat n_s$
 and angle $2\pi b_n t$.
  Then $\xi^0(t)[z_n:\,z_{n+2}]=[w_t]$ and
$\xi^1(t)[z_n:\,z_{n+2}]=[u_t]$.
  We write
 $[w_{s,t}]:=\xi^s(t)[z_n:\,z_{n+2} ].$

 For each $s\in[0,1]$ we define
 $$\psi^s(t)[z]=\big([w_{s,t}],\, [(w_{s,t})^{a_j}z_j\epsilon_j:\,z_{n+1}]
 \big).$$
Thus $\psi^0=\psi_{\bf b}$ and  $\psi^1=\psi_{\bf b'}$. Since
$\xi^s(t+t')=\xi^s(t)\circ\xi^s(t')$, the family $\{\psi^s\}$ is a
homotopy  between $\psi_{\bf b}$ and $\psi_{\bf b'}$ consisting of
circle actions.
 \qed

\subsection{ Hirzebruch surfaces.} \label{Section4}

\medskip

Here we will consider the manifold ${\mathbb
P}(L_1\oplus\dots\oplus L_{n-1}\oplus{\mathbb C })$ when $n=2$. In
this case $\Delta$ is the trapezoid in $({\mathbb R}^2)^*$ with
vertices $P_1(0,0)$, $P_2(0,\tau)$, $P_3(\sigma,0)$,
$P_4(\sigma,\tau+a\sigma)$ and $M_{\Delta}$ is a Hirzebruch
surface ${\mathbb P}(L_k\oplus {\mathbb C})$, where $L_k$ the
holomorphic line bundle over ${\mathbb C}P^1$ with first Chern
number $k:=-a$. More explicitly
 $$M_{\Delta}=\{(z_j)\in{\mathbb C}^4\,: \,
|z_1|^2+|z_3|^2=\sigma/\pi,\;
-a|z_1|^2+|z_2|^2+|z_4|^2=\tau/\pi\}/\simeq,$$ with
$$(z_1,\,z_2,\,z_3,\,z_4)\simeq(\epsilon_1\epsilon_2^{-a}z_1,\,
\epsilon_2z_2,\, \epsilon_1z_3,\, \epsilon_2z_4),$$
 $\epsilon_i$ being an element of $U(1)$.


 The map $[z_j]\mapsto ([z_2:z_4],[z_2^{-a}z_3:z_4^{-a}z_3:z_1])$ gives a
 representation of $M$ as a submanifold of ${\mathbb C}P^1\times {\mathbb
 C}P^2$.

An element ${\bf b}$ in the lattice ${\mathbb Z}^2$ determines a
Hamiltonian circle action $\psi_{\bf b}$ on $M_{\Delta}$.
The following Proposition is a particular case of Theorem
\ref{TmequivaPL}

\begin{Prop}\label{Propgh(0a)}
 If ${\bf b}=(b_1,\,b_2)$ in the integer lattice  ${\mathbb
 Z}^2$ and ${\bf b'}=(b_1-kb_2,\, -b_2)$,
 then $\psi_{\bf b}\sim \psi_{\bf b'}$.
 \end{Prop}

\smallskip

Now {\em we assume that
 $k=-a\ne 0$}. The case $k=0$ will be considered in Subsection
 \ref{Sectspheres}. In order to calculate $S_{\bf b}:=S_{\psi_{\bf b}}$ it is
 necessary to take into account the Type of ${\bf b}$.
  The element ${\bf b}\ne 0$ is of Type(0) when $\langle l,\,{\bf b}\rangle\ne 0$, for all edge $l$ of $\Delta$. And
  it is of Type (1) if there is an edge $f$ such that $\langle f,\,{\bf
  b}\rangle=0$.

\smallskip

{\it Values of $S_{-}$ on ${\rm Type}(0)$}

 If ${\bf b}=(b_1,\,b_2)\in{\mathbb Z}^2$ is of type $0$, then $b_1\ne 0\ne b_2\ne 0\ne
b_0:=(b_1+ab_2)$, and the fixed points of the $S^1$-action
$\psi_{\bf b}$ are the inverse images by the moment map $\mu$ of
the vertices of $\Delta$. If $P$ is a vertex, the
$U(1)$-equivariant Euler class of the normal bundle of
$\mu^{-1}(P)$ in $M$ is the product of the weights of the
corresponding isotropy representation. On the other hand, the
edges meeting at $P$ can be expressed as $\{P+t\rho_i\,:\, t\in
[0,c_i] \}$, where $\rho_i\in({\mathbb Z}^2)^*$ and $\rho_1,
\rho_2$ is a basis of $({\mathbb Z}^2)^*$. The aforementioned
weights are precisely the numbers $\langle \rho_i,\, {\bf
b}\rangle$, with $i=1,2$.

The normalized Hamiltonian function for the $S^1$-action
$\psi_{\bf b}$ is $\langle\mu,\, {\bf b}\rangle-\langle{\rm
Cm}(\Delta),\,{\bf b}\rangle$, where ${\rm Cm}(\Delta)$ is the
center of mass of $\Delta$. Hence
\begin{equation}\label{Sbfb}
S_{\bf b}=\frac{e^{q u}}{u^2}\,\sum_{j=1}^{4}\frac{e^{-\langle
P_j,\,{\bf b}\rangle u}}{\prod_{i=1}^2\langle \rho_{ji},\,{\bf
b}\rangle},
\end{equation}
with $\{\rho_{ji}\}_{i=1,2}$ is the basis of $({\mathbb Z}^2)^*$
formed with the aforementioned vectors which give the directions
of the edges meeting at $P_j$, and $q:= \langle{\rm
Cm}(\Delta),\,{\bf b}\rangle.$

On $U_1:=\{[z]\in M\,:\, z_3\ne 0\ne z_4  \}$ we can consider the
complex coordinates $w=z_2/z_4$ and $w'=z_1/(z_3z_4^k)$. Then
$$\psi_{\bf b}(w,w')=(e^{2\pi ib_2t}w,\, e^{2\pi
ib_1t}w').$$
 The point $\mu^{-1}(P_1)$ has coordinates $w=0=w'$, so the
isotropy representation at this point has weights $b_1$ and $b_2$;
  in fact these values are  the products by ${\bf b}$ of the direction
vectors of edges meeting at $P_1$. The remaining terms in
(\ref{Sbfb}) can be obtained in a similar way, or directly from
the vectors $\rho_i$. So (\ref{Sbfb}) can be written
\begin{equation}\label{sbfHz1}
u^2S_{\bf b}= e^{qu} \big(A- B e^{-\tau b_2u}- A
 e^{-\sigma b_1u} + B  e^{-(\tau b_2+\sigma
b_0)u}  \big)
\end{equation}
where $A:=(b_1b_2)^{-1}$ and $B:=(b_0b_2)^{-1}$. Note that $A\ne
B$; otherwise $b_0$ would be equal to $b_1$, but this is
impossible in ${\rm Type}(0)$.

The expression (\ref{sbfHz1}) will be written
\begin{equation}\label{sbfHz2}
u^2S_{\bf b}=  Ae^{\alpha_1 u}- Be^{\alpha_2 u}- Ae^{\alpha_3 u}
  + B e^{\alpha_4 u},
\end{equation}
where $\alpha_1=q,\; \alpha_2=q-\tau b_2,\; \alpha_3=q-\sigma
b_1,\; \alpha_4=q-(\tau b_2+\sigma b_0).$ If $\alpha_2\ne\alpha_3$
and $\alpha_1\ne \alpha_4$, there are four different exponential
functions in (\ref{sbfHz2}). By contrast, if $\alpha_2=\alpha_3$
and $\alpha_1=\alpha_4$ there are only two different exponentials
in (\ref{sbfHz2}). Therefore we will distinguish three subtypes in
${\rm Type}(0)$.
$$
{\rm Type}(0) \begin{cases} {\bf b}\in{\rm Subtype}(0)\alpha\;\; {\rm iff}\,\, \alpha_2\ne\alpha_3,\alpha_1\ne\alpha_4  \\
{\bf b}\in{\rm Subtype}(0)\beta\;\; {\rm iff}\;\,{\rm either}\;
\alpha_2=\alpha_3,\alpha_1\ne\alpha_4,\;\,
{\rm or}\; \alpha_2\ne\alpha_3,\alpha_1=\alpha_4  \\
{\bf b}\in{\rm Subtype}(0)\gamma\;\; {\rm iff}\;\,
\alpha_2=\alpha_3,\alpha_1=\alpha_4
 \end{cases}
$$

It follows from the Lemma \ref{Lemexponbisbis}  the following
Proposition

\begin{Prop}\label{Propsubtype1}
If ${\bf b}$ and ${\bf b'}$ are in ${\rm Type}(0)$ and belong to
different subtype, then $S_{\bf b}\ne S_{\bf b'}$.
\end{Prop}

 Given ${\bf b}, {\bf b'}\in{\rm Subtype}(0)\alpha$, if $S_{\bf
 b}=S_{\bf b'}$,
then by Lemma \ref{Lemexponbisbis} it follows from (\ref{sbfHz2})
 \begin{equation}\label{EqualitAB}
 \{ Ae^{\alpha_1u},\,-Be^{\alpha_2u},\,-Ae^{\alpha_3u},\,Be^{\alpha_4u}\}=
 \{ A'e^{\alpha'_1u},\,-B'e^{\alpha'_2u},\,-A'e^{\alpha'_3u},\,B'e^{\alpha'_4u}\}.
 \end{equation}
 This set equality gives, in principle,
the following eight possibilities

$${\bf b},{\bf b'}\in{\rm Subtype (0)}\,\alpha\begin{cases}

  (i)\; A=A' \; \;\;\;\;\, \begin{cases}(a)\; \; B=B'\\
                        (b)\;\; B=-B'
           \end{cases} \\
\smallskip

(ii) \;A=B' \;\;\;\;\, \begin{cases}(a)\;\; B=A'\\
                        (b)\;\; B=-A'
           \end{cases}\\
\smallskip

(iii) \; A=-A' \; \begin{cases}(a)\;\; B=B'\\
                        (b)\;\; B=-B'
           \end{cases}\\

\smallskip

(iv) \; A=-B' \;\,  \begin{cases}(a)\;\; B=A'\\
                        (b)\; \;B=-A'
           \end{cases}

\end{cases}$$
Where   the corresponding exponential function which multiplies
each constant  has been deleted in the  equalities; for example,
we have written   $B=-A'$ instead of  $Be^{\alpha_4
u}=-A'e^{\alpha'_3u}$.

\begin{Prop}\label{Subtyp0alpha}
If ${\bf b}$ and ${\bf b'}$ are in ${\rm Subtype}(0)\alpha$ and
$S_{\bf b}=S_{\bf b'}$, then   ${\bf b}\equiv {\bf b'}.$
 \end{Prop}

{\it Proof.} In the case (i)(a),
$\,A=A',\,B=B',\,\alpha_1=\alpha'_1,\,\alpha_4=\alpha'_4.$ And by
(\ref{EqualitAB}) $-Be^{\alpha_2 u}=-B'e^{\alpha'_2 u},$
$\,-Ae^{\alpha_3 u}=-A'e^{\alpha'_3 u}$. From $\alpha_j=\alpha'_j$
for all $j$, we deduce ${\bf b'}={\bf b}.$

In the case (i)(b), $-Ae^{\alpha_3 u}=-A'e^{\alpha'_3 u}$. So
$\alpha_3=\alpha'_3$ and of course $\alpha_1=\alpha'_1$. Thus
$b_1=b'_1$. From $A=A'$ it follows $b'_2=b_2$; that is, ${\bf
b}={\bf b'}$. But this is contradictory with $B=-B'$.

It is also easy to check that the cases (ii)(a), (ii)(b),
(iii)(a), (iii)(b) and (iv)(a) do not occur. In the case (iv)(b)
it is straightforward to deduce $b'_1=b_1-kb_2$ and $b'_2=-b_2.$
\qed

 \smallskip

Next we study the map $S_{-}$ on ${\rm Subtype}(0)\beta$. If ${\bf
b}$ belongs to ${\rm Subtype}(0)\beta$, then either
 $$ u^2S_{\bf b}=
Ae^{\alpha_1 u} - (A+B)e^{\alpha_2 u}
  + B e^{\alpha_4 u},\;\;{\rm when}\;
   \alpha_2=\alpha_3$$
   or
$$ u^2S_{\bf b}=
(A+B)e^{\alpha_1 u} - Ae^{\alpha_2 u}-Be^{\alpha_3 u}
   ,\;\;{\rm when}\;
   \alpha_1=\alpha_4.$$

We will consider the partition ${\mathcal P}$ of ${\rm
Subtype}(0)\beta$ formed by two subsets: The one  consisting of
all elements for which $\alpha_2=\alpha_3$, and the other one its
complement subset. It is not hard to prove the following
Proposition
 \begin{Prop}\label{PropSub(0)beta}
 Given ${\bf b},{\bf b'}\in {\rm Subtype}(0)\beta$ with $S_{\bf b}=S_{\bf b'}$, if they belong
 to the same subset of the partition ${\mathcal P}$, then ${\bf b'}={\bf b}.$
 If they belong to different subsets, then ${\bf b'}\equiv{\bf b}.$
 \end{Prop}

It is easy to prove
\begin{Prop}\label{PropSub(0)gamma}
 The map $S_{-}$ is injective on   ${\rm Subtype}(0)\gamma$.
 \end{Prop}

Propositions  \ref{Subtyp0alpha}, \ref{PropSub(0)beta} and
\ref{PropSub(0)gamma} can be put together

\begin{Prop}\label{PropType0}
If ${\bf b},{\bf b'}\in{\rm Type}(0)$ and $S_{\bf b}=S_{\bf b'}$,
then ${\bf b}\equiv{\bf b'}.$
 \end{Prop}

\smallskip

{\it Values of $S_{-}$ on ${\rm Type}(1)$}

If ${\bf b}$ belongs to  ${\rm Type}(1)$
  one of $b_j$'s vanishes. We will
distinguish three possibilities: $b_2=0$, $\,b_1=0$ and ${\bf
b}=(kb_2,\,b_2)$, that is $b_0=0$.

We first consider the case $b_2=0$. Now $\psi_{\bf b}[z]=[e^{2\pi
ib_1t}z_1,z_2,z_3,z_4].$ So the fixed point set is $\{[z]\in
M\,:\,z_1z_3=0\}$. On $U_2=\{[z]\in M\,:\, z_2\ne 0\ne z_3\}$ we
consider the complex coordinates $v=z_4/z_2$ and
$v'=z_1/(z_3z_2^k)$. So $\frac{\partial}{\partial v'}$ is a
section of the normal bundle $N_Q$ to $Q=\{[z]: z_1=0\}$ on $Q\cap
U_2$. On the other hand,  $\frac{\partial}{\partial w'}$ is a
section of $N_Q$ on $Q\cap U_1$. As
$$\frac{\partial}{\partial w'}=\Big(\frac{z_4}{z_2}\Big)^k\frac{\partial}{\partial
v'},$$
  identifying $Q$ with ${\mathbb C}P^1$ it turns out that $N_Q$ is
${\mathcal O}(k).$

 Since the coordinate $w'$ is transformed by
$\psi_{\bf b}(t)$ in $\,{\rm exp}(2\pi ib_1t)w',\,$ the
$U(1)$-equivariant Euler class $E(N_Q)$ of the normal bundle to
$N_Q$ is $c_1({\mathcal O}(k))+b_1u$.

Analogously, if we denote by $Q':=\{[z]: z_3=0\}$, then
$E(N_{Q'})=c_1({\mathcal O}(-k))-b_1u$. As $\omega(Q)=\tau$ and
$\omega(Q')= \tau-k\sigma=:\lambda,$
 it follows from the localization formula
 \begin{equation}\label{sbfb(2)a}
S_{\bf b}=e^{qu}\Big( \frac{1}{b_1 u}\big(\tau-\frac{k}{b_1 u}
\big)- \frac{e^{-\sigma b_1 u}}{b_1 u}\big(\lambda-\frac{k}{b_1 u}
\big) \Big),
\end{equation}
where $q:=\langle{\rm Cm}(\Delta),\,{\bf b}\rangle$.

\smallskip

Next we consider the case    when ${\bf b}=(0,\,b_2)$, with
$b_2\ne 0$. Now the fixed point set is
$\{[z]\,:\,z_2=0\}\cup\{P_2,P_4\}$. The vector field
$\frac{\partial}{\partial w}$ is a global section to the normal
bundle to $\{[z]\,:\, z_2=0\}$, so its $U(1)$-equivariant Euler
class is $b_2u$. The equivariant Euler classes of the normal
bundles $N_{P_2}$ and $N_{P_4}$ are $k(b_2 u)^2$ and $-k(b_2 u)^2$
respectively. As $\omega(\{[z]\,:\, z_2=0\})=\sigma$, it follows
\begin{equation}\label{sbfb(2)b1}
S_{\bf b}=e^{qu}\Big( \frac{e^{-\tau b_2 u}}{k(b_2u)^2 }
-\frac{e^{-\lambda b_2 u}}{k(b_2u)^2 } + \frac{\sigma}{b_2u}
  \Big).
\end{equation}

\smallskip

Finally, let us suppose that ${\bf b}=(kb_2,\,b_2)$. The
corresponding fixed point set is
$\{[z]\,:\,z_4=0\}\cup\{P_1,P_3\}$. The equivariant Euler class of
the normal bundles to $\{[z]\,:\,z_4=0\}$, $P_1$ and $P_3$ are
$-b_2u$, $k(b_2u)^2$ and $-k(b_2u)^2$ respectively. So

\begin{equation}\label{sbfb(2)b2}
S_{\bf b}=e^{qu}\Big( \frac{1}{k(b_2u)^2} -\frac{e^{-\sigma k b_2
u}}{k(b_2u)^2 } - \frac{\sigma e^{-\tau  b_2 u}}{b_2u}
  \Big).
\end{equation}

To compare   the expressions  $S_{\bf b}$ for different elements
${\bf b}$,
we classify the elements of ${\rm Type}(0)$
 according to the following scheme
of subtypes
$$
{\rm Type}(1)\;\begin{cases} {\bf b}\in{\rm Subtype (1)}\alpha \;\; \text{iff}\;\;b_2=0 \\
                    {\bf b}\in {\rm Subtype (1)}\beta\;\;\text{iff}\;\; b_2\ne 0
         \end{cases}
 $$

 From  (\ref{sbfb(2)a}), (\ref{sbfb(2)b1})
 and (\ref{sbfb(2)b2}) together with Lemma \ref{Lemexponbisbis} it
 follows the following Proposition
\begin{Prop}\label{Propsubtype}
 If   ${\bf b},$ ${\bf b'}$
 belong to ${\rm Type} (1)$ but to different subtype, then $S_{\bf b}\ne S_{\bf b'}$.
\end{Prop}

\smallskip

It follows from Theorem \ref{ThmHHA1} together with Theorem
\ref{ThmGEnralG}, Proposition \ref{Propsubtype} and Proposition
\ref{Propsubtype1} the following Theorem

\begin{Thm}\label{Thmimtro2}
Let $M_{\Delta}$ be a Hirzebruch surface non diffeomorphic to
$S^2\times S^2$. Let ${\bf b}$ and ${\bf b'}$ be  nonzero elements
of ${\mathbb Z}^2$. If they belong to different subtype then
$\psi_{\bf b}\nsim\psi_{\bf b'}.$
\end{Thm}

Next we study whether the map ${\bf b}\mapsto S_{\bf b}$ is
injective restricted to each of the different subtypes of ${\rm
Type}(1).$

If ${\bf b}=(b_1,\,0)$ belongs to ${\rm Subtype}(1)\alpha$, in the
expression of $u^2S_{\bf b}$ obtained from (\ref{sbfb(2)a}) there
are two different exponentials
$$u^2S_{\bf b}=A(u)e^{qu}+B(u)e^{(q-\sigma b_1)u},$$
where $A(u)$ and $B(u)$ are the following polynomials
$A(u)=(b_1)^{-1}\tau u-k(b_1)^{-2},\,$ $B(u)=-(b_1)^{-1}\lambda
u+k(b_1)^{-2}.$ If ${\bf b'}$ belongs to Subtype (1)$\alpha$ and
$S_{\bf b}=S_{\bf b'}$ then either $A(u)=A'(u)$ and $B(u)=B'(u)$
or $A=B'$ and $B=A'$. In the first case $b_1'=b_1$. In the second
one implies $-(b_1)^2=(b'_1)^2.$ Hence we have the following
Proposition
\begin{Prop}\label{PropinyecS1a}
$S_{-}$ is injective on Subtype(1)$\alpha$.
\end{Prop}

In ${\rm Subtype}(1)\beta$ one can distinguish the case in which
${\bf b}$ is of the form $(0,c)$ and when it is of the form $(kc,
c)$. As in the preceding Proposition, if ${\bf b}=(0,\,b_2)$,
$\,{\bf b'}=(0,\,b'_2)$ and $S_{\bf b}=S_{\bf b'}$, then ${\bf
b}={\bf b'}$. Analogously  when ${\bf b}=(kc,\,c)$ and ${\bf
b'}=(kc',\,c')$.

Finally, given ${\bf b}=(0,\,b_2)$ and ${\bf
b'}=(b'_1,\,b'_2)=(kc,\, c)$, in the expressions for $u^2S_{\bf
b}$ and $u^2S_{\bf b'}$ obtained from (\ref{sbfb(2)b1}) and
(\ref{sbfb(2)b2}) there are three different exponentials
$$u^2S_{\bf b}=Ae^{(q-\tau b_2)u}+B e^{(q-\lambda
b_2)u}+C(u)e^{qu},$$
$$u^2S_{\bf b'}=A'e^{q'u}+B' e^{(q'-\sigma
k b'_2)u}+C'(u)e^{(q'-\tau b'_2)u},$$
 with $C(u)=b_2^{-1}\sigma$ and $C'(u)=-(b'_2)^{-1}\sigma$.
 The equality $S_{\bf b}=S_{\bf b'}$ implies $C(u)=C'(u)$; that
 is $b'_2=-b_2$. In other words ${\bf b}\equiv{\bf b'}.$
  We have
 proved the Proposition

  \begin{Prop}\label{PropSubt1beta}
If ${\bf b}$ and ${\bf b'}$ belong to ${\rm Subtype}(1)\beta$ with
$S_{\bf b}=S_{\bf b'}$ then
  ${\bf b}\equiv{\bf b'}$.
  \end{Prop}

The preceding Propositions can be summarize in the following one

\begin{Prop}\label{PropType1}
If ${\bf b}$ and ${\bf b'}$ are  elements   of  type $1$ with
$S_{\bf b}=S_{\bf b'}$ then ${\bf b}\equiv{\bf b'}.$
\end{Prop}

{\bf Proof of Theorem \ref{Thmmain}.}   Theorem \ref{Thmmain}  is
consequence of Proposition \ref{PropType1}, Proposition
\ref{PropType0}, Proposition \ref{Propgh(0a)} and Theorem
\ref{ThmHHA1}.
 \qed

\subsection{The manifold $S^2\times S^2$.}\label{Sectspheres}

Here we study the relation $\psi_{\bf b}\sim\psi_{\bf b'}$
  when $M$ is the symplectic
manifold associated to the polytope $\Delta$ in the case $a=0$.
That is, $M$  is diffeomorphic to $S^2\times S^2$. Now $b_0=b_1$
  and ${\rm Cm}(\Delta)=(\sigma/2,\,\tau/2)$, so the
corresponding discussion is simpler than the one of case $a\ne 0$.

If $b_1\ne 0\ne b_2$, the expression (\ref{sbfHz1}) (which is also
valid when $k=0$) reduces to
\begin{equation}\label{Sbsp1}
S_{\bf b}=\frac{4}{b_1b_2u^2}\Big(\sinh \big( \frac{\sigma b_1
u}{2} \big) \sinh\big( \frac{\tau b_2 u}{2} \big)  \Big)
\end{equation}

When ${\bf b}=(b_1,0)$ the expression for $S_{\bf b}$ can be
obtained from (\ref{sbfb(2)a})
\begin{equation}\label{Sbsp12}
S_{\bf b}=\frac{2\tau}{b_1u} \sinh \big( \frac{\sigma b_1 u}{2}
\big)
\end{equation}
Analogously, if ${\bf b}=(0,b_2)$, then
\begin{equation}\label{Sbsp3}
S_{\bf b}=\frac{2\sigma}{b_2u} \sinh \big( \frac{\tau b_2 u}{2}
\big)
\end{equation}

Using Lemma \ref{Lemexponbisbis} it is straightforward to deduce
from (\ref{Sbsp1}) the following Proposition
\begin{Prop}\label{Prok01}
If $b_1b_2\ne 0\ne b'_1b'_2$, then $S_{\bf b}=S_{\bf b'}$ iff one
of the following statements are true:

 (i) $\,|b'_1|=|b_1|$ and $|b'_2|=|b_2|$.

 (ii) $\;\sigma b_1=\pm\tau b'_2$ and $\tau b_2=\pm\sigma b'_1.$

 (iii) $\;\sigma b_1=\pm\tau b'_2$ and $\tau b_2=\mp\sigma b'_1.$
  \end{Prop}

{\bf Proof of Theorem \ref{Thmspheres}.}     In this
case$$\psi_{\bf b}(t)[z]=\big([e^{2\pi
ib_2t}z_2\,:\,z_4],\,[z_3\,:\,e^{2\pi ib_1t}z_1 ]\big),$$ where
$[z]=\big([z_2:z_4],\,  [z_3:z_1]\big)\in S^2\times S^2$. On the
other hand, if $\sigma/\tau\notin{\mathbb Q}$, the possibilities
(ii) and (iii) in the preceding Proposition do not occur. If, for
example $b'_1=-b_1$ and $b'_2=b_2$, let $\xi^s_t$ be the rotation
of angle $2\pi b_1t$ around the vector $\Hat n_s$ introduced in
the Proof of   Theorem \ref{TmequivaPL}. We set
$$\psi^s_t[z]= \big([e^{2\pi
ib_2t}z_2\,:\,z_4],\,\xi^s_t[z_3\,:\, z_1 ]\big).$$
 Then $\psi^s$ is a homotopy between $\psi_{\bf b}$ and $\psi_{\bf
 b'}$ consisting of $U(1)$-actions. That is, $|b'_1|=|b_1|$ and
 $|b'_2|=|b_2|$ implies $\psi_{\bf b}\sim\psi_{\bf b'}$.
So   Theorem \ref{Thmspheres} follows from Proposition
\ref{Prok01} together with Theorem \ref{ThmHHA1}.
 \qed

\smallskip

{\it Remark.}
 In general $S_{\bf b}=S_{\bf b'}$ does not imply
$[\psi_{\bf b}]=[\psi_{\bf b}]\in[U(1),\,{\rm Ham}(M)]$. For
instance, if $\tau=\sigma$, ${\bf b}=(1,0)$ and ${\bf b'}=(0,1)$,
it follows from (\ref{Sbsp12}) and (\ref{Sbsp3}) that $S_{\bf
b}=S_{\bf b'}$. Gromov proved that ${\rm Ham}(S^2\times S^2)$
deformation retracts to the group $SO(3)\times SO(3)$ \cite{Gr}.
So $\pi_1({\rm Ham}(S^2\times S^2))={\mathbb Z}_2\times {\mathbb
Z}_2$ and $\psi_{\bf b}$ and  $\psi_{\bf b'}$ are the generators
of the respective groups in the product. Hence $\psi_{\bf b}$ and
$\psi_{\bf b'}$ are not homotopic, and a fortiori $\psi_{\bf
b}\nsim \psi_{\bf b'}$.


\section{$S^1$-actions in coadjoint orbits.}\label{Coadjoint}

Let $G$ be a compact connected Lie group and $\eta\in{\frak g}^*$.
We consider the coadjoint orbit ${\mathcal O}$ of the element
$\eta$, endowed with the standard symplectic structure $\omega$
(see \cite{aK76}). A moment map $\mu_G:{\mathcal O}\to {\frak
g}^*$ for the action of $G$ on ${\mathcal O}$ is  the opposite of
the inclusion.
  The manifold ${\mathcal O}$ can be
identified with the quotient $G/G_{\eta}$, where $G_{\eta}$ is the
subgroup of isotropy of $\eta$. Let $H$ be a Cartan subgroup of
$G$ contained in $G_{\eta}$ (\cite{vG84}, page 166).
We have the decomposition of ${\frak g}_{\Bbb C}$ as direct sum of
root spaces
$${\frak g}_{\Bbb C}={\frak h}_{\Bbb
C}\oplus\bigoplus_{\alpha\in{\mathcal R}}{\frak g}_{\alpha},$$
where ${\mathcal R}$ is the set of roots determined by $H$. We
will denote by $\check\alpha\in{\frak h}$ the coroot associated to
$\alpha\in{\mathcal R}$ \cite{G-W}. Let ${\frak p}$ be the
parabolic subalgebra
$${\frak p}={\frak h}_{\Bbb C}\oplus\bigoplus_{\eta(i\check
\alpha)\leq 0}{\frak g}_{\alpha}.$$
 If $P$ is the parabolic
subgroup of $G_{\Bbb C}$ generated by
 ${\frak p}$, then $G_{\Bbb C}/P$ is a complexification of
 ${\mathcal O}$ and
 \begin{equation} \label{T10}
 T^{1,0}_{\eta}{\mathcal O}=\bigoplus_{\alpha\in{\mathcal
 R}^+}{\frak g}_{\alpha},
 \end{equation}
where ${\mathcal R}^+=\{\alpha\in{\mathcal R}\,:\,\langle \eta,
i\check\alpha\rangle>0 \}$.

Given $X\in{\frak g}$ we consider the map $h:=\langle
\mu_G,\,X\rangle:{\mathcal O} \to {\Bbb R},$ and define the
constant $\kappa\in{\Bbb R}$ by the relation
$$\kappa\int_{\mathcal O}\omega^n=\int_{\mathcal O} h\,\omega^n.$$
The function $h-\kappa$ is the normalized Hamiltonian function for
the isotopy in ${\mathcal O}$ generated by $X$ through the
$G$-action.

If $X$ belongs to
  the integer lattice of ${\frak
h}$, then $X$ determines a $U(1)$-action $\psi$ on ${\mathcal O}$.
By $u$ we denote a coordinate in ${\frak u}(1)$, then closed
$U(1)$-equivariant $2$-form $\omega-(\mu-\kappa)u$ is a
representative of the coupling class $c_{\psi}$ \cite{K-M}. Thus
$$S_{\psi}=e^{\kappa u}\int_{\mathcal O}e^{-\mu
u+\omega}\in\Hat H(BU(1)).$$

In order to apply the localization formula we need to determine
the fixed point set of the $S^1$-action $\psi$. Let us assume that
$X$ is  a regular element \cite{wF91} of
  the lattice of ${\frak h}$. We denote by $W_{\eta}$ the stabilizer subgroup of $\eta$ in
the Weyl group $W$ of $(G,H)$. It is easy to see that the fixed
point set of the action $\psi$ is $\{w\cdot\eta \,:\, w\in
W/W_{\eta} \}$ (see \cite{B-G-V} page 231).

Given $w\in W$, the $U(1)$-equivariant Euler class of the normal
bundle to the point $w\cdot\eta$ in ${\mathcal O}$ can be
calculated as follows:
 From  (\ref{T10}) it follows
 $$T^{1,0}_{w\cdot\eta}{\mathcal O}=\bigoplus_{\alpha\in{\mathcal
R}^+}{\frak g}_{w\cdot\alpha}.$$
   Moreover $[Y,\,Z]=(w\cdot\alpha)(Y)Z$ for any $Z\in{\frak
 g}_{w\cdot\alpha}$ and all $Y\in{\frak h}$.
  So the representation of $U(1)$ on ${\frak g}_{w\cdot\alpha}$
induced by $\psi$ therefore has the  weight
$-i(w\cdot\alpha)\in{\frak u}(1)^*$, and the weights  of the
isotropy representation at the point $w\cdot \eta\in{\mathcal O}$
are
 $-i(w\cdot\alpha)$, with $\alpha\in{\mathcal R}^+$. That is,   the
 $U(1)$-equivariant Euler class satisfies
 \begin{equation}\label{eqEu2}
 E(N_{w\cdot\eta})(Y)=\prod_{\alpha\in{\mathcal R}^+}-i(w\cdot\alpha)(Y),
 \end{equation}
 for any $Y=tX$ with $t\in{\Bbb R}$. Thus $E(N_{w\cdot\eta})$ is a
 monomial in $u$
 \begin{equation}\label{Eulcoa}
 E(N_{w\cdot\eta})=\prod_{\alpha\in{\mathcal R}^+}A^{w}_{\alpha}u\in
 H(BU(1)),
 \end{equation}
  with $A^w_{\alpha}$ constant.

   Consequently, if $X$ is a regular element of the integer lattice of
   ${\frak h}$,  in the localization
  $\Hat H(BU(1))_u$ one has
  \begin{equation}\label{Sphi}
  S_{\psi}=e^{\kappa u}\sum_{w\in
  W/W_{\eta}}\frac{e^{\eta(w^{-1}\cdot X)u}}{\prod_{\alpha\in{\mathcal R}^+} A^w_\alpha
  u},
  \end{equation}
since $\mu(w\cdot\eta)=-\langle
w\cdot\eta,\,X\rangle=-\eta(w^{-1}\cdot X)$.

\begin{Thm}\label{TheoremCoadjo}
Given $X,X'$ regular elements of the integer lattice of ${\frak
h}$, let $\psi$ and $\psi'$ be the circle actions on ${\mathcal
O}$ generated by $X$ and $X'$, respectively. If there is no a
translation in ${\Bbb R}$ which applies the set of real numbers
$\{\eta(w\cdot X)\,:\, w\in W/W_{\eta}\}$ in $\{\eta(w\cdot
X')\,:\, w\in W/W_{\eta}\}$, then $\psi$ and $\psi'$ are not
homotopically equivalent by means of a homotopy consisting of
circle actions
\end{Thm}

{\it Proof.}  The class $S_{\psi'}$ is obtained from (\ref{Sphi})
by substituting $\kappa$, $X$ and $A^{w}_{\alpha}$ by the
corresponding $\kappa'$, $X'$ and $A'^{w}_{\alpha}$. By the
hypothesis the sets of functions
$${\mathcal E}:=\{e^{(\kappa +\eta(w^{-1}\cdot X))u}\,|\,  w\in W/W_{\eta}\}\;\;{\rm and}\;\;
{\mathcal E'}:=\{e^{(\kappa' +\eta(w^{-1}\cdot X'))u}\,|\,  w\in
W/W_{\eta}\} $$ are distinct. Hence ${\mathcal E}$ is not
contained in the ${\Bbb R}[u]$-submodule of ${\Bbb R}[[u]]$
generated by ${\mathcal E'}$. Reasoning as in  the proof of
 Theorem
\ref{Thmtorus}   we conclude that
 $S_{\psi}\ne S_{\psi'}$, and the theorem
follows from Theorem \ref{ThmHHA1}.

\qed

\begin{Prop}\label{Propbeta}
Let $\psi$ and $\psi'$ be the circle actions on the coadjoint
orbit ${\mathcal O}$ generated by the regular elements $X$ and
$X'$ of the integer lattice of ${\frak h}$, respectively.
 If there is no an affine transformation
 $$B:x\in{\mathbb R}\to\lambda x+a\in{\mathbb
 R},\;\;{\rm with}\; 0\ne\lambda\in{\mathbb
 Z}\;{\rm and}\;a\in{\mathbb R},$$
such that
$$B\big(\{\eta(w\cdot X)\,:\, w\in W/W_{\eta}\}\big)=
\{\eta(w\cdot X')\,:\, w\in W/W_{\eta}\},$$
 then $\psi$ and
$\psi'$ are $\sim_{rp}$-inequivalent.
\end{Prop}

{\it Proof.} If $\psi$ and $\psi'$ were $\sim_{rp}$-equivalent,
there would be an automorphism $v$ of $U(1)$ with
$\psi'\sim\psi\circ v$. Moreover $\psi\circ v$ is the $S^1$-action
generated by $Z=c X$, with $c$ a nonzero integer. By Theorem
\ref{TheoremCoadjo} there must be a translation $T_b$ in ${\mathbb
R}$ which applies $\{\eta(w\cdot Z)= c\,\eta(w\cdot X)\,:\, w\in
W/W_{\eta}\}$ in $\{\eta(w\cdot X')\,:\, w\in W/W_{\eta}\}$. Hence
the affine transformation $B:x\in{\mathbb R}\mapsto c
x+b\in{\mathbb R}$ contradicts the hypothesis of the Proposition.

 \qed

\smallskip

{\bf Proof of Theorem \ref{ThmCoadj}.}
 We consider the coadjoint orbit ${\mathcal O}=G/Q$ of $G$, where
$Q$ is a subgroup of $G$ which contains a maximal torus $H$ of
$G$. The element $X $ in the integer lattice of ${\frak h}$
  generates the circle action
  $$\psi_t:gQ\in G/Q\mapsto (e^{tX}g)Q\in G/Q.$$
We denote by $\psi'$ the $S^1$-action defined by   $X'$.
 Let $g_1$ be an element of the normalizer $N(H)$ of $H$ in
  $G$, such that $g_1H\in N(H)/H=W$ corresponds to the
  element $w$. And let $g: s\in[0,1]\mapsto g(s)\in
   G$ be a curve with $g(0)=1$ and $g(1)=g_1$.
   We
  write $X^s:={\rm Ad}g(s)X,$  then
  $X^1=w(X)=X'$.
For $s, t\in[0,\,1]$ we put
$$\psi^s_t:gQ\in G/Q\mapsto \big(e^{tX^s}g\big)Q\in G/Q.$$
As $e^{X^s}=g(s)e^X(g(s))^{-1}=1$, for each $s$ the family
$\{\psi^s_t\}_t$ is a circle action on the orbit ${\mathcal
O}=G/Q$. Moreover
$$\psi^1_t(gQ)=\big(e^{tX^1} g \big)Q =\big(e^{tX'} g \big)Q=\psi'_t(gQ),$$
and $\psi_t^0(gQ)=\psi_t(gQ)$. Hence $\{\psi^s\}$ is a homotopy
 consisting of $U(1)$-actions between $\psi$ and $\psi'$.

\qed

\smallskip
By $G_k({\mathbb C}^n)$ we denote the Grassmann manifold
consisting of the $k$-subspaces in ${\mathbb C}^n$. If $M$ is a
coadjoint orbit of $SU(n)$ diffeomorphic to $G_k({\mathbb C}^n)$,
it is the orbit of an element $\eta\in{\frak su}(n)^*$ such that
for $C\in{\frak su}(n)$
\begin{equation}\label{Grassmann}
\eta (C)=\alpha_1\sum_{j=1}^k  C_{jj} + \alpha_2\sum_{j=k+1}^n
C_{jj},
 \end{equation}
  with $\alpha_1,\alpha_2\in i{\mathbb R}$ and $\alpha_1\ne
\alpha_2$. That is, $\eta(C)=-ir\sum_{j=1}^k C_{jj},$ with
$r\in{\mathbb R}\setminus \{0\}.$

\smallskip

Now we suppose that  $k=1$, then
  the coadjoint orbit ${\mathcal
O}:={\mathcal
 O}_{\eta}=SU(n)/U(n-1)={\mathbb C}P^{n-1}$.   The coset $gU(n-1)\in
 SU(n)/U(n-1)$ determines the line in ${\mathbb
 C}^{n}$  defined by the first column of $g$.
  This is
 the identification ${\mathcal O}\simeq {\mathbb
 C}P^{n-1}$ we will use.

 We write ${\mathcal R}$ for the set of roots
determined by $H$, the diagonal subgroup of $SU(n)$. The
stabilizer of $\eta$ in the Weyl group $W=N(H)/H$  consists of all
the permutations $\sigma$ of $\{1,\dots, n  \}$ with
$\sigma(1)=1$. Denoting by $\sigma_1$ the
 permutation identity and by $\sigma_j$ the transposition $(1,j)$,
with $j=2,\dots, n$, then $\{\sigma_1,\dots,\sigma_n\}$ is a set
of representatives of the classes in $W/W_{\eta}$.

Given $X\in{\frak su}(n)$, the Hamiltonian isotopy generated by
$X$ has the map $h-\kappa$ as normalized Hamiltonian function,
$\kappa$ being a constant and $h$
 \begin{equation}\label{Hamilnonorm}
 h:\nu\in{\mathcal O}\mapsto
-\langle\nu,\,X\rangle\in{\mathbb R}.
 \end{equation}

If $X=i(a_1,\dots,a_n)$ with $(a_j)\in (2\pi{\mathbb Z})^{n}$,
then $X$ generates a circle action $\psi$ on ${\mathcal O}$.
Furthermore, if $a_i\ne a_j$ for $i\ne j$, then the fixed point
set for the action $\psi$ is
$\{\sigma_k\cdot\eta\}_{k=1,\dots,n}$.

Denoting by $\alpha_{ij}$ the root defined by
$\alpha_{ij}(B)=b_i-b_j$, for $B={\rm diag}\,(b_k)$, then the set
of positive roots is
 ${\mathcal R}^+=\{\alpha_{12},\,\dots,\alpha_{1n}\}$,
assumed $r>0$. For $k\ne 1$
$$(-i\sigma_k\cdot\alpha_{1s})(X)=\begin{cases}
a_k-a_s,\;{\rm if}\; s\ne k \\
a_k-a_1 ,\;{\rm if}\; s= k\end{cases}. $$ And
$(-i\sigma_1\cdot\alpha_{1s})(X)=a_1-a_s$. By (\ref{eqEu2}) we
obtain for the equivariant Euler class of normal bundle $N_k$ to
$\sigma_k\cdot\eta$
  \begin{equation}\label{Eul2}
   E(N_k)=\prod_{s\ne k}(a_k-a_s) u^{n-1}.
   \end{equation}

   Since
$h(\sigma_k\cdot \eta)=-\eta(\sigma_k^{-1}\cdot X)=-ra_k$, it
follows from (\ref{Eul2}) together with (\ref{Localiz})
\begin{equation}\label{Spsi}
S_{\psi}=\frac{e^{\kappa
u}}{u^{n-1}}\sum_{k=1}^n\frac{e^{ra_ku}}{\prod_{s\ne
k}(a_k-a_s)}\;.
 \end{equation}

\begin{Prop}\label{Thcp} Let $X=(x_1,\dots,x_n),\,X'=(x'_1,\dots,x'_n)$ be regular elements of the integer
lattice of ${\frak h}$, and $\psi$, $\psi'$ respective the circle
actions on the coadjoint orbit ${\mathcal O}$ of $SU(n)$
diffeomorphic to ${\mathbb C}P^{n-1}$. Then the following
statements are equivalent:

 (a) $S_{\psi}=S_{\psi'}.$

(b) There is a permutation $\tau$ of $\{1,\dots,n \}$ such that
$x'_j=x_{\tau(j)}$, for all $j$.

(c) $\psi\sim\psi'$.

\end{Prop}

{\it Proof.}  (a) $\,\Longrightarrow\,$ (b). If we write
$X'=i(a_1',\dots, a'_n)$ the expression for $S_{\psi'}$ can be
obtained from (\ref{Spsi}) by substituting the $a_j$'s by the
$a'_j$'s and $\kappa$ by $\kappa'$. From the equality
$S_{\psi}=S_{\psi'}$ together with the ${\mathbb
R}[u,u^{-1}]$-linear independence of the exponentials, it follows
that there is a permutation $\tau$ of $(1,\dots,n)$ such that
$\kappa'+ra'_j=\kappa+ra_{\tau(j)}$. Thus there is a constant $c$
such that $a'_j-a_{\tau(j)}=c$ for all $j=1,\dots, n$. As $\sum _j
a'_j=\sum _i a_i=0$, it turns out $c=0$; hence (b) holds.

(b) $\,\Longrightarrow\,$ (c) by  Corollary \ref{CoadjointSUn}.
 The equivalence of (a), (b) and (c)
follows from Theorem \ref{ThmHHA1}.
 \qed

\smallskip

 Theorem \ref{ThmCPn} is the generalization of Proposition \ref{Thcp} when
 the hypothesis about the regularity of $X$ and $X'$ is deleted.

{\bf Proof of Theorem \ref{ThmCPn}.}

If
\begin{equation}\label{EqiX}
-iX= \big( \overbrace{a,\dots, a}^{m}, \overbrace{b,\dots,
b},\dots ,\overbrace{d,\dots, d}
\big)\equiv(a_1,\dots,a_n)\in(2\pi{\mathbb Z})^{n}
\end{equation}
with $a,b,\dots$ different numbers, then a connected component of
the fixed point set for the circle action $\psi$ determined by $X$
is $F=\{[z_1:\dots:z_{m}:0:\dots:0]\}\subset{\mathbb C}P^{n-1}=
{\mathcal O}.$ The contribution of the component $F$ to $S_{\psi}$
can be calculated as in the paragraphs  before Lemma
\ref{Lemgene}. That is, $N_F$ can be decomposed in a direct sum of
$U(1)$-equivariant line bundles $\oplus_{j=m+1}^n N_j$, and
$$E(N_F)=\prod_{j=m+1}^n\big( c_1(N_j)+(a-a_j)u  \big).$$
So we have as in (\ref{Eqparadespues})
$$(E(N_F))^{-1}=\frac{1}{u^{n-m}}\big(\alpha_0+\frac{\alpha_1}{u}+\dots
\big).$$
 The aforementioned contribution to $S_{\psi}$  is
 $$e^{(\kappa+ra)u}Q_a(u),$$
$Q_a(u)$ being a Laurent polynomial of degree $m-n$. Thus
expression for $S_{\psi}$ is
\begin{equation}\label{Eqpsix}
 S_{\psi}=e^{(\kappa+ra)u}Q_a(u)+\dots +e^{(\kappa+rd)u}Q_d(u).
\end{equation}

Let
 \begin{equation}\label{Eqix'}
  -iX'= \big( \overbrace{a',\dots, a'}^{m'},  \dots ,\overbrace{e',\dots, e'}
\big)\equiv(a'_1,\dots,a'_n)\in(2\pi{\mathbb Z})^{n}.
 \end{equation}
 The equality $S_{\psi}=S_{\psi'}$ implies the following set
 equality $\{a+\beta,b+\beta,\dots, d+\beta  \}=\{a',\dots, e' \}$,
 with $\beta$  constant. Furthermore, if for example $a'=b+\beta$
 then
 ${\rm degree}\,Q_{a'}= {\rm degree}\,Q_{b}.$ That is, the number
 of $b$'s in (\ref{EqiX}) must be equal to the number of $a'$'s in
(\ref{Eqix'}).  Hence there is a permutation $\tau$ of
$\{1,\dots,n\}$,  such that $a_j+\beta=a'_{\tau(j)}$ and if
$a_i=a_j$ then $a'_{\tau(i)}=a'_{\tau(j)}.$ As
$\sum_ja_j=0=\sum_ia'_i$, then $\beta=0$. So the theorem is
consequence of Corollary \ref{CoadjointSUn} and Theorem
\ref{ThmHHA1}.
 \qed

\smallskip

{\bf Proof of Theorem \ref{ThmGrass}.}

 Now we must consider the case  $1<k<n$ in (\ref{Grassmann}). In this case, if
 the permutation $\sigma\in W_{\eta}$, then $\sigma$ applies the
 subset $\{1,\dots,k\}$ onto itself (and of course $\{k+1,\dots,n\}$ onto itself).
 The stabilizer $W_{\eta}$ has $k!(n-k)!$ elements. A set of
 representatives of the elements of $W/W_{\eta}$ is the set ${\mathcal
 C}_{k}$ of  combinations   of $\{1,\dots,n  \}$ with $n-k$
 elements. If ${\frak d}\in {\mathcal C}_k$, $\sigma_{\frak d}$ will
 denote the corresponding element in the quotient $W/W_{\eta}$

Let $X={\rm diag}(ia_1,\dots,ia_n)\in (2\pi i{\mathbb Z})^n$ be a
 regular element, and $\psi $ is the corresponding circle action on the
Grassmannian $G_k({\mathbb C}^n)$. If ${\frak c}=\{1,\dots,
n\}\setminus\{j_1,\dots,j_k\}$,  the value of the  Hamiltonian
(\ref{Hamilnonorm}) on ${\frak c}\cdot\eta$ is
$$h(\sigma_{\frak c}\cdot\eta)=-r\sum_{l=1}^ka_{j_l}.$$
We write $\alpha_{\frak d}=\sum_{l=1}^ka_{i_l}$, when ${\frak
d}=\{i_1,\dots,i_k\}\in{\mathcal C}_k$. Then
\begin{equation}\label{SGrass}
 S_{\psi}=\frac{e^{\kappa u}}{u^{k(n-k)}}\sum_{{\frak c}\in{\mathcal C}_k}
\frac{e^{r\alpha_{\frak c}u }}{p_{\frak c}},
 \end{equation}
where $p_{\frak c}\in{\mathbb Z}$ and $r$ is a real number.

If $X'$ is another regular element and $S_{\psi}=S_{\psi'}$, it
follows from (\ref{SGrass}) that there is
 a bijective map $f:{\mathcal C}_k\to {\mathcal C}_k$ such that
$$\kappa+\alpha_{f(\frak c)}=\kappa'+\alpha'_{\frak c},$$
for all  ${\frak c}\in{\mathcal C}_k$.
\qed

\smallskip

If ${\mathcal O}$ is a coadjoint orbit of $SU(n)$ diffeomorphic to
the full flag manifold, it is the orbit of an element
$\eta\in{\frak su}(n)^*$, such that for any $C\in{\frak su}(n)$
\begin{equation}\label{etafull}
\eta(C)=-i\sum_{j=1}^{n-1}r_jC_{jj},
\end{equation}
with $r_i\ne r_j$ for $i\ne j$. We will assume that
$r_1>r_2>\dots>r_{n-1}.$ So ${\mathcal R}^+=\{\alpha_{ij}\,:\,
i<j\}.$

Given $X\in{\frak su}(n)$ with $-iX=(a_1,\dots,a_n)\in\big(2\pi
{\mathbb Z}\big)^n$ and $a_i\ne a_j$ for $i\ne j$, the fixed point
set for the $S^1$-action $\psi$ defined by $X$ on ${\mathcal O}$
is $\{\sigma\cdot\eta\,:\,\sigma\in{\mathcal S}_n \}$, where
${\mathcal S}_n$ is the symmetric group. The $U(1)$-equivariant
Euler class
$$E(N_{\sigma\cdot\eta})=u^m\epsilon(\sigma)\prod_{r<s}(a_r-a_s),$$
$\epsilon(\sigma)$ being the signature of $\sigma$, and
$m=n(n-1)/2$. As
$$-\eta(\sigma\cdot X)=\sum_{j=1}^{n-1}r_ja_{\sigma(j)},$$
we obtain
\begin{equation}\label{Spsifullflag}
S_{\psi}=\frac{e^{\kappa
u}}{u^m\prod_{r<s}(a_r-a_s)}\sum_{\sigma\in{\mathcal
S}_n}\epsilon(\sigma)\, {\rm
exp}\Big(u\sum_{j=1}^{n-1}r_ja_{\sigma(j)}  \Big).
\end{equation}


If $X'$ is another regular element with $-iX'=(a'_1,\dots,
a_n')\in (2\pi{\mathbb Z})^n$ and $S_{\psi}=S_{\psi'}$, it follows
from (\ref{Spsifullflag}) together with Lemma \ref{Lemexponbisbis}
that there exist a constant $\beta$ and a bijective map $f$ from
${\mathcal S}_n$ into itself, such that for all
$\sigma\in{\mathcal S}_n$
\begin{equation}\label{Eqfullflag}
\sum_{j=1}^{n-1}r_ja'_{\sigma(j)}=\sum_{j=1}^{n-1}r_ja_{f\sigma(j)}+\beta.
\end{equation}
We have proved the following Proposition

\begin{Prop}\label{Fullflag}
Let ${\mathcal O}$ be the coadjoint orbit of the element $\eta\in
{\frak su}(n)^*$ given by (\ref{etafull}), and let $\psi$ and
$\psi'$ be the $U(1)$-actions defined by the regular elements
$(ia_1,\dots,ia_n)$ and $(ia'_1,\dots,ia'_n)$ respectively. If
$\psi\sim\psi'$, then there exist a constant $\beta$ and a
bijective map $f$ from ${\mathcal S}_n$ into itself, such that
(\ref{Eqfullflag}) holds for all $\sigma\in{\mathcal S}_n$.
\end{Prop}

{\it Remark.} In particular, if there exist $\tau\in{\mathcal
S}_n$, such that $a'_j=a_{\tau(j)}$ for all $j$; that is, the
hypothesis of Corollary \ref{CoadjointSUn} are satisfied, then we
can construct the map $f:\sigma\in{\mathcal S}_n\mapsto
\tau\cdot\sigma\in {\mathcal S}_n.$ Thus
$a'_{\sigma(k)}=a_{f\sigma(k)}$ and equation (\ref{Eqfullflag})
holds with $\beta=0$. That is, Proposition \ref{Fullflag} is
consistent with Corollary \ref{CoadjointSUn}.




\medskip


\end{document}